# Metrics of positive Ricci curvature on quotient spaces

Lorenz J. Schwachhöfer[*] and Wilderich Tuschmann[†]

## 1 Introduction

One of the classical problems in differential geometry is the investigation of closed manifolds which admit Riemannian metrics with given lower bounds for the sectional or the Ricci curvature and the study of relations between the existence of such metrics and the topology and geometry of the underlying manifold. Despite many efforts during the past decades, this problem is still far from being understood. For example, so far the only obstructions to the existence of a metric with positive Ricci curvature come from the obstructions to the existence of metrics with positive scalar curvature ([Li], [Hi], [Ro], [SchY], [Ta]) and the Bonnet-Myers theorem which implies that the fundamental group of a closed manifold with positive Ricci curvature must be finite.

Fruitful constructions of metrics with positive Ricci curvature on closed manifolds have so far been established by techniques that include deformation of metrics ([Au], [Eh], [We]), Kähler geometry ([Yau1], [Yau2]), bundles and warping ([Po], [Na], [BB], [GPT]), special kinds of surgery ([SY], [Wr]), metrical glueing ([GZ2]), and Sasakian geometry ([BGN]). Particularly large classes of examples of manifolds with positive Ricci curvature are given by all compact homogeneous spaces with finite fundamental group ([Na]) and all closed cohomogeneity one manifolds with finite fundamental group ([GZ2]).

In this article we present several new classes of manifolds which admit metrics of positive or nonnegative Ricci curvature. The idea is to consider quotients of manifolds $(M, g)$ of positive or nonnegative Ricci curvature by a free isometric action. While taking such a quotient non-decreases the sectional curvature, it may well happen in general that the Ricci curvature of the quotient is inferior to the one of $M$. Our first results, however, state that the quotient of $M$ *does* admit metrics of positive Ricci curvature if $M$ belongs to one of the aforementioned classes of homogeneous spaces or spaces with a cohomogeneity one action. More precisely, we prove:

**Theorem A** *A biquotient $G/\!/H$ of a compact Lie group $G$ admits a Riemannian metric of positive Ricci curvature which is invariant under the canonical action of the normalizer $Norm_{G\times G}(H)/H$ if and only if its fundamental group is finite.*

**Theorem B** *Let $(M, G)$ be a closed cohomogeneity one manifold, and suppose that $L \subset G$ acts freely on $M$. Then the quotient $M/L$ admits a Riemannian metric of positive Ricci curvature which is invariant under the canonical action of the normalizer $Norm_G(L)/L$ if and only if its fundamental group is finite.*

We shall say that a closed manifold $M$ admits metrics of *positive (nonnegative, respectively) Ricci and almost nonnegative sectional curvature* if for every $\varepsilon > 0$ there is a metric $g_\varepsilon$ on $M$ such that $Ric(M, g_\varepsilon) > 0$ ($Ric(M, g_\varepsilon) \geq 0$, respectively) and $Sec(M, g_\varepsilon) \cdot diam(M, g_\varepsilon)^2 > -\varepsilon$.

Admitting metrics of positive Ricci and almost nonnegative sectional curvature is a much stronger property than the mere existence of a metric with positive Ricci curvature. Indeed, by [Gr] and [SY] there are in each dimension $n \geq 4$ infinitely many closed manifolds with positive Ricci curvature which do not admit almost nonnegatively curved metrics. Also, for closed manifolds of positive Ricci and almost nonnegative sectional curvature there is no known obstruction to the existence of metrics of nonnegative sectional curvature.

[*]Research supported by the Communauté française de Belgique, through an Action de Recherche Concertée de la Direction de la Recherche Scientifique.

[†]Research supported by a DFG Heisenberg Fellowship.

As we shall show, the biquotients in Theorem A actually admit metrics of positive Ricci and almost nonnegative sectional curvature. This might appear to be a minor point since biquotients obviously carry metrics of nonnegative sectional curvature; it is unclear, however, if every biquotient with finite fundamental group admits a metric of positive Ricci *and* nonnegative sectional curvature.

Moreover, by [ST], any closed cohomogeneity one manifold and hence any of the quotients considered in Theorem B admits also invariant metrics of almost nonnegative sectional curvature. However, it is not always clear if positive Ricci and almost nonnegative sectional curvature can be achieved simultaneously, though there are large classes of quotients where we can show that this is the case (cf. Corollary 6.2 and Remark 6.7). In particular, this holds when $L = 1$, so that we obtain as a slight improvement of [GZ2] and [ST]:

**Theorem C** *Let $(M,G)$ be a closed cohomogeneity one manifold. Then $M$ admits $G$-invariant metrics of positive Ricci and almost nonnegative sectional curvature if and only if its fundamental group is finite.*

Theorem $C$ is in fact optimal within the class of invariant metrics, since there are closed simply connected cohomogeneity one manifolds which do not admit invariant metrics of nonnegative sectional curvature ([GVWZ]).

Let us now drop the condition of finiteness of the fundamental group. Any closed cohomogeneity one manifold admits invariant metrics of nonnegative Ricci curvature ([GZ2]) and metrics of almost nonnegative sectional curvature ([ST]). We generalize these results as follows.

**Theorem D** *Let $(M,G)$ be a closed cohomogeneity one manifold, and suppose that $L \subset G$ acts freely on $M$. Then there exist metrics of nonnegative Ricci and almost nonnegative sectional curvature on the quotient $M/L$ which are invariant under the canonical action of the normalizer $Norm_G(L)/L$.*

Theorems $A$, $B$ and $C$ yield many new examples of manifolds with positive Ricci and almost nonnegative sectional curvature. For instance, we have the following classes.

**Theorem E** *For every integer $m \geq 2$, there are families $(N_d^{4m-2})_{d \geq 1}$ and $(\tilde{N}_d^{8m-4})_{d \geq 1}$ of mutually not homotopy equivalent closed simply connected manifolds of dimension $4m-2$ and $8m-4$, respectively, admitting metrics of positive Ricci and almost nonnegative sectional curvature which are invariant under an action of cohomogeneity two. Their integral cohomology groups and cohomology rings are given by*

$$H^k\left(N_d^{4m-2}\right) \cong H^k(\mathbb{CP}^{2m-1}) \qquad and \qquad H^k\left(\tilde{N}_d^{8m-4}\right) \cong H^k(\mathbb{HP}^{2m-1}) \quad for\ all\ k, \quad and$$

$$H^*(N_d^{4m-2}) \cong H^*(\tilde{N}_d^{8m-4}) \cong \mathbb{Z}[x,y]/\{x^m - d\ y,\ y^2\}.$$

*Thus, these manifolds have the rational cohomology ring of $\mathbb{CP}^{2m-1}$ and $\mathbb{HP}^{2m-1}$, respectively. Moreover, for $d \geq 3$ these manifolds neither admit a Lie group action of cohomogeneity less than two, nor are they diffeomorphic to a biquotient. Finally, there are fiber bundles $S^2 \hookrightarrow N_d^{8m-2} \to \tilde{N}_d^{8m-4}$.*

For $m = 2$, the first class was already considered in ([GZ1], Corollary 3.9) where it was shown that the manifolds $N_d^6$ are $S^2$-bundles over $S^4$ and admit metrics of even nonnegative sectional curvature.

The examples in Theorem E are obtained by considering quotients of free $S^1$- and $S^3$-actions on certain *Brieskorn manifolds*. The latter establish a class of odd dimensional cohomogeneity one manifolds which includes all Kervaire spheres and infinitely many rational homology spheres in each dimension $n = 4m - 1$, $m \geq 2$. By Theorem C, they all admit invariant metrics of positive Ricci and almost nonnegative sectional curvature. There are also interesting quotients of the Brieskorn manifolds by finite cyclic groups which are again of cohomogeneity one ([ST]).

**Theorem F** *There are at least $4^k$ oriented diffeomorphism types of homotopy $\mathbb{RP}^{4k+1}$ which admit a cohomogeneity one action and thus invariant metrics of positive Ricci and almost nonnegative sectional curvature.*

*For every integer $m \geq 3$, the Kervaire spheres admit free actions of $\mathbb{Z}_m$ whose quotient is again of cohomogeneity one and thus admits invariant metrics of positive Ricci and almost nonnegative sectional curvature.*

An independent construction of Sasakian metrics of positive Ricci curvature on the above mentioned homotopy $\mathbb{RP}^{4k+1}$'s has been established recently in [BGN]. Moreover, the four homotopy $\mathbb{RP}^5$'s even admit invariant metrics of nonnegative sectional curvature by [GZ1]. Notice also that the second family in Theorem $F$ yields homotopy lens spaces which are differentiably distinct from the standard ones in those dimensions where the Kervaire spheres are exotic.

The remaining parts of this paper are structured as follows. In section 2, we consider biquotients and prove the result which implies Theorem $A$. In section 3, we set up some notation and prove a deformation result which is needed later on. Then we recall the so-called Cheeger trick in section 4 which we apply to obtain a straightforward proof of the existence of metrics with almost nonnegative sectional curvature on cohomogeneity one manifolds which is due to B.Wilking. We also obtain further results which are of importance later. In sections 5 and 6 we work towards the proof of Theorems $B$, $C$ and $D$. Finally, in section 7 we consider the Brieskorn manifolds and their quotients and, in particular, the concrete examples which yield Theorems $E$ and $F$.

It is our pleasure to thank Karsten Grove, Burkhard Wilking and Wolfgang Ziller for helpful discussions and suggestions, and the Max-Planck-Institute for Mathematics in the Sciences in Leipzig for its hospitality.

## 2 Ricci curvature of biquotients

Let $G$ be a compact Lie group with Lie algebra $\mathfrak{g}$, and let $H \subset G \times G$ be a closed subgroup for which the action of $H$ on $G$ given by

$$(h_1, h_2) \cdot g := h_1 g h_2^{-1}.$$

is free and hence the quotient $G/\!/H$ is a manifold, called the *biquotient of $G$ by $H$*.

Recall that for groups $G_2 \subset G_1$, the *normalizer of $G_2$ in $G_1$* is defined as

$$Norm(G_2) = Norm_{G_1}(G_2) := \{g_1 \in G_1 \mid Ad_{g_1}(G_2) = G_2\},$$

where we omit the subscript if no confusion is possible. Note that $G_2 \subset Norm_{G_1}(G_2)$ is a normal subgroup.

Evidently, conjugation yields a canonical action of $Norm_{G \times G}(H)$ on $G/\!/H$, and $H$ acts trivially by this action. Thus, there is a canonical induced action of $Norm(H)/H$ on $G/\!/H$.

The tangent spaces of the fibers of the projection $\pi : G \to G/\!/H$ are given by

$$\mathcal{V}_g = \left(Ad_{g^{-1}} \circ pr_1 - pr_2\right)(\mathfrak{h}),$$

where $pr_i : \mathfrak{g} \oplus \mathfrak{g} \to \mathfrak{g}$ is the projection onto the $i$th factor and where we identify $T_g G \cong \mathfrak{g}$ via left translation. We fix a biinvariant inner product $Q$ on $\mathfrak{g}$ and denote the $Q$-orthogonal complement of $\mathcal{V}_g$ by

$$\mathcal{H}_g := (\mathcal{V}_g)_Q^\perp.$$

Since the action of $H$ on $G$ preserves $Q$, there is a unique Riemannian metric $g_Q$ on $G/\!/H$ for which the projection $\pi : G \to G/\!/H$ becomes a Riemannian submersion. Evidently, such a metric has always nonnegative sectional curvature and is invariant under the action of $Norm(H)/H$.

**Theorem 2.1** Let $M := G/\!/H$ be a biquotient of a compact Lie group $G$.

1. There is a fibration $M \to T^k$ whose fiber is a biquotient with finite fundamental group, where $T^k$ is a torus of dimension $k = \dim(\mathfrak{z}(\mathfrak{g}) \cap \mathcal{H}_g)$, $g \in G_0$ and where $\mathfrak{z}(\mathfrak{g})$ denotes the center. In particular, this dimension is independent of $g \in G_0$, and $k = 0$ iff $M$ has finite fundamental group.

2. The following are equivalent.

   (a) $M$ admits $Norm(H)/H$-invariant Riemannian metrics of positive Ricci and almost nonnegative sectional curvature.

   (b) $M$ has finite fundamental group.

   (c) For any biinvariant metric $Q$ on $\mathfrak{g}$, the submersion metric $(M, g_Q)$ has positive Ricci curvature on a dense open subset of $M$.

**Proof.** To show the first part, we may assume that $G$ is connected and – after replacing $G$ by a finite cover if necessary – that $G = G_s \times Z$ where $G_s$ is semisimple and $Z$ is a torus whose Lie algebra is $\mathfrak{z}(\mathfrak{g})$. The fibers of the principal $H$-bundle $\pi : G \to G/\!/H$ are the images of the maps

$$\imath^g : H \longrightarrow G, \qquad \imath^g(h_1, h_2) := h_1 g h_2^{-1}.$$

We also consider the maps

$$\jmath^g : H \longrightarrow Z, \qquad \jmath^g := pr_Z \circ L_{g^{-1}} \circ \imath^g,$$

where $pr_Z : G = G_s \times Z \to Z$ is the projection. Thus,

$$\jmath^g(h_1, h_2) = pr_Z(h_1) pr_Z(h_2)^{-1}$$

is independent of $g \in G$ and a Lie group homomorphism.

Let $G' := G_s \times \jmath^g(H)$. Then evidently, $G'$ is a normal subgroup such that $G/G' = T^k$ is a torus of dimension $k = \text{codim}(\jmath^g(H) \subset Z) = \text{codim}(d\jmath^g(\mathfrak{h}) \subset \mathfrak{z}(\mathfrak{g})) = \text{codim}(pr_{\mathfrak{z}(\mathfrak{g})}(\mathcal{V}_g) \subset \mathfrak{z}(\mathfrak{g})) = \dim(\mathfrak{z}(\mathfrak{g}) \cap \mathcal{H}_g)$. Note that $\jmath^g(H)$ and hence $k$ is independent of $g \in G$. Moreover, we may assume that $H \subset G' \times G'$ so that the canonical projection $G \to G/G' = T^k$ induces a fibration $M = G/\!/H \to T^k$ whose fiber is $G'/\!/H$.

Now the map $\jmath^g_* : \pi_1(H) \to \pi_1(\jmath^g(H))$ has finite cokernel, being induced by a Lie group epimorphism. Since $L_{g^{-1}}$ is a diffeomorphism and $\pi_1(G_s)$ is finite, it follows that the map $\imath^g_* : \pi_1(H) \to \pi_1(G')$ has finite cokernel as well. Now the homotopy exact sequence of the principal $H$-bundle $G' \to G'/\!/H$ implies that $\pi_1(G'/\!/H)$ is finite, which shows the first part.

For the second part, we do no longer assume $G$ to be connected.
(a) $\Rightarrow$ (b) follows from the Bonnet-Myers theorem.
(c) $\Rightarrow$ (a) follows from the deformation results in [Au], [Eh] and [We], and the fact that for any metric on $M = G/\!/H$ which is induced by a biinvariant metric on $G$, $Norm(H)/H$ acts by isometries.
(b) $\Rightarrow$ (c) For $g \in G$, we let

$$\mathcal{F}_g := \{v \in \mathcal{H}_g \mid [v, \mathcal{H}_g] = 0\} \subset \mathcal{H}_g.$$

The dimension of $\mathcal{F}_g$ is semicontinuous, whence there is an element $g_0 \in G$ around which this dimension is constant. W.l.o.g. we may assume that $g_0 = e$ after conjugating one of the components of $H$. We let

$$S := Stab(\mathcal{F}_e) = \{g \in G \mid Ad_{g|\mathcal{F}_e} = Id_{\mathcal{F}_e}\} \subset G$$

be the stabilizer of $\mathcal{F}_e \subset \mathfrak{g}$. Then $S$ is compact, and the Lie algebra $\mathfrak{s}$ of $S$ is the centralizer of $\mathcal{F}_e$, i.e.

$$\mathfrak{s} = \{v \in \mathfrak{g} \mid [v, \mathcal{F}_e] = 0\}.$$

Thus, $\mathcal{H}_e \subset \mathfrak{s}$. Moreover, we let

$$H_S := H \cap (S_0 \times S_0),$$

where $S_0$ is the identity component of $S$. Now there is a smooth map $p : S_0/\!/H_S \to G/\!/H$ which makes the following diagram commute, where the map on the top is the inclusion:

$$\begin{array}{rccc} \imath: & S_0 & \longrightarrow & G \\ & \downarrow \pi_S & & \downarrow \pi \\ p: & S_0/\!/H_S & \longrightarrow & G/\!/H \end{array} \qquad (1)$$

Our goal shall be to show that $p$ is a finite covering map, and we now finish the proof in several steps.

*Step 1 :* $\mathcal{F}_e \subset \mathcal{H}_g$ for all $g \in S$.

To see this, we let $h \in H$, $g \in S$ and $a \in \mathcal{F}_e$. Then $Ad_g a = a$ by definition of $S$, and we calculate

$$\begin{aligned} Q(a, (Ad_{g^{-1}} \circ pr_1 - pr_2)(h)) &= Q(Ad_g a, pr_1(h)) - Q(a, pr_2(h)) \\ &= Q(a, (pr_1 - pr_2)(h)) = 0, \end{aligned}$$

where the first equality follows from the biinvariance of $Q$ and the last from $(pr_1 - pr_2)(h) \in \mathcal{V}_e$ and $a \in \mathcal{F}_e \subset \mathcal{H}_e$.

*Step 2 :* $\mathcal{F}_g \subset \mathfrak{s}$ for all $g \in S$.

This follows immediately since by step 1, $[\mathcal{F}_g, \mathcal{F}_e] \subset [\mathcal{F}_g, \mathcal{H}_g] = 0$.

*Step 3 :* $\mathcal{H}_g \subset \mathfrak{s}$ for all $g \in S_0$.

First of all, by the local constancy of the rank $\mathcal{F}_g \to \mathcal{F}_e$ as subspaces as $g \to e$, and $[\mathcal{F}_e, x] \neq 0$ for all $0 \neq x \in \mathfrak{s}^\perp$. Thus, we have $[\mathcal{F}_g, x] \neq 0$ for all $0 \neq x \in \mathfrak{s}^\perp$ and all $g \in S$ in a sufficiently small neighborhood of the identity.

If $x_\mathfrak{s} + x_{\mathfrak{s}^\perp} \in \mathcal{H}_g$ and $a \in \mathcal{F}_g$, then $0 = [a, x_\mathfrak{s}] + [a, x_{\mathfrak{s}^\perp}]$ and since $a \in \mathfrak{s}$ by step 2, the first summand lies in $\mathfrak{s}$ while the second lies in $\mathfrak{s}^\perp$, i.e. both summands vanish simultaneously. Thus $[\mathcal{F}_g, x_{\mathfrak{s}^\perp}] = 0$, which by the preceding remark implies that $x_{\mathfrak{s}^\perp} = 0$ and hence $\mathcal{H}_g \subset \mathfrak{s}$ for all $g \in S$ sufficiently close to the identity.

By analyticity, this implies that $\mathcal{H}_g \subset \mathfrak{s}$ for all $g \in S_0$.

*Step 4 :* The differentials $dp$ are surjective.

Since $T_{\pi(g)} G/\!/H = d\pi(\mathcal{H}_g)$ and $\mathcal{H}_g \subset \mathfrak{s}$ for all $g \in S_0$ by step 3, it follows that $d(\pi \circ \imath) : \mathfrak{s} \cong T_g S_0 \to T_{\pi_g} G/\!/H$ is surjective, and by the commutativity of (1), the assertion follows.

*Step 5 :* The differentials $dp$ are injective.

By step 4, $\dim \ker dp$ is constant, hence it suffices to verify the injectivity of $dp$ at $\pi_S(e)$. Since $\pi_S$ is a submersion, it suffices to show that $\ker(d(p \circ \pi_S)) = \ker(d\pi_S)$, and by the commutativity of (1) this is equivalent to

$$\mathfrak{s} \cap (pr_1 - pr_2)(\mathfrak{h}) = (pr_1 - pr_2)(\mathfrak{h} \cap (\mathfrak{s} \oplus \mathfrak{s})). \qquad (2)$$

The inclusion $\supset$ is obvious. Suppose that $(h_1, h_2) \in \mathfrak{h}$ with $h_1 - h_2 \in \mathfrak{s}$. Then for all $a \in \mathcal{F}_e$ and $(x_1, x_2) \in \mathfrak{h}$, we have

$$\begin{aligned} Q([a, h_2], x_2) &= Q(a, [h_2, x_2]) \quad \text{since } Q \text{ is biinvariant} \\ &= Q(a, [h_1, x_1]) \quad \text{since } [h_1, x_1] - [h_2, x_2] \in \mathcal{V}_e \text{ and } a \in \mathcal{F}_e \subset \mathcal{H}_e \\ &= Q([a, h_1], x_1) \quad \text{since } Q \text{ is biinvariant} \\ &= Q([a, h_2], x_1) \quad \text{since } [a, h_1 - h_2] \in [\mathcal{F}_e, \mathfrak{s}] = 0. \end{aligned}$$

Therefore, $Q([a, h_2], x_1 - x_2) = 0$ for all $(x_1, x_2) \in \mathfrak{h}$, i.e. $[a, h_2] \in \mathcal{H}_e$, hence $[a, [a, h_2]] \in [\mathcal{F}_e, \mathcal{H}_e] = 0$. Thus, $0 = Q([a, [a, h_2]], h_2) = -Q([a, h_2], [a, h_2])$ again by the biinvariance of $Q$, whence $[a, h_2] = 0$ for all $a \in \mathcal{F}_e$, i.e. $h_2 \in \mathfrak{s}$, and hence $(h_1, h_2) \in \mathfrak{h} \cap (\mathfrak{s} \oplus \mathfrak{s})$ which proves (2).

*Step 6 :* $p$ is a finite covering map.

By steps 4 and 5 it follows that $p$ is a local diffeomorphism, and $S_0$ and thus $S_0/\!/H_S$ are compact connected.

*Step 7 :* There is a dense open subset $U \subset G$ such that $\mathcal{F}_g = 0$ for all $g \in U$.

By hypothesis (b), $G/\!/H$ has finite fundamental group, and hence by step 6 so does $S_0/\!/H_S$. Thus, by the first part of the theorem, this implies that $\mathfrak{z}(\mathfrak{s}) \cap \mathcal{H}_e = 0$, and since by definition of $S$ this space contains $\mathcal{F}_e$, it follows that $\mathcal{F}_e = 0$.

On the other hand, by the very first remark of the proof (b) $\Rightarrow$ (c), we can replace $e$ by any $g_0 \in G$ around which $\dim \mathcal{F}_g$ is constant, hence the claim follows.

*Step 8 :* $\pi(U) \subset G/\!/H$ is a dense open subset on which the Ricci curvature of $g_Q$ is positive. Thus, (c) holds.

$\pi$ is a submersion, hence $\pi(U) \subset G/\!/H$ is dense open as $U \subset G$ is dense open. Since $(M, g_Q)$ has nonnegative sectional curvature, we have $Ric(v_p) = 0$ for $v_p \in T_pM$ iff $Sec(v_p \wedge w_p) = 0$ for all $w_p \in T_pM$. Identifying $\mathcal{H}_g \cong T_{\pi(g)}M$ via the differential $d\pi$, O'Neill's formula implies that $[v_g, \mathcal{H}_g] = 0$, where $v_g \in \mathcal{H}_g$ is such that $d\pi(v_g) = v_p$, i.e. $v_g \in \mathcal{F}_g$. If $g \in U$ then this implies that $v_g = 0$ and thus $v_p = 0$, i.e. on $\pi(U)$ the Ricci curvature of $g_Q$ is positive. ∎

It is not clear whether or not the biinvariant metric on a closed biquotient with finite fundamental group itself has positive Ricci curvature. We do not know of any example where this is *not* the case, but an affirmative answer is known only in very special cases ([Fa]).

## 3 $G$-invariant metrics

Let $H \subset G$ be compact Lie groups, and let $\mathfrak{h} \subset \mathfrak{g}$ be their Lie algebras. Fix a biinvariant inner product $Q$ on $\mathfrak{g}$ which yields the $Q$-orthogonal decomposition

$$\mathfrak{g} = \mathfrak{h} \oplus \mathfrak{m}.$$

Let $\varphi : \mathfrak{m} \to \mathfrak{m}$ be an $Ad(H)$-equivariant linear map which is symmetric and positive definite w.r.t. $Q$, i.e. such that the bilinear form on $\mathfrak{m}$ given by

$$g_\varphi(x, y) := Q(x, \varphi y)$$

yields an inner product on $\mathfrak{m}$. Clearly, we can extend $\varphi$ to $\mathfrak{g}$ by setting $\varphi|_\mathfrak{h} = Id_\mathfrak{h}$ and thus obtain an $Ad(H)$-equivariant inner product on $\mathfrak{g}$ which in turn induces a left invariant, $Ad(H)$-invariant Riemannian metric on $G$, also denoted by $g_\varphi$. Also, there is a unique $G$-invariant Riemannian metric on the homogeneous space $G/H$ such that the natural projection $\pi : (G, g_\varphi) \to G/H$ becomes a Riemannian submersion. By abuse of notation, we denote this metric by $g_\varphi$ as well.

This procedure establishes a one-to-one correspondence between $G$-invariant Riemannian metrics on $G/H$ and symmetric, positive definite $Ad(H)$-equivariant linear maps $\varphi : \mathfrak{m} \to \mathfrak{m}$. Moreover, we call the metrics induced by $\varphi = c_0 Id_\mathfrak{m}$ with $c_0 > 0$ *normal homogeneous*.

A $G$-invariant metric on $M := I \times G/H$ is determined by a smooth one-parameter family $\varphi(t) : \mathfrak{m} \to \mathfrak{m}$, $t \in I$, such that the metric has the form

$$g = g_\varphi = dt^2 + g_{\varphi(t)}, \tag{3}$$

where the metric $g_{\varphi(t)}$ on $G/H$ is induced by $\varphi(t)$. Then the curvature of such a $g$ has been calculated in [GZ2] and [ST] as follows.

**Proposition 3.1** *Let $M = I \times G/H$ and $g = dt^2 + g_{\varphi(t)}$ be as in (3), and let $c \in \mathbb{R}$, $x, y \in T_{eH}G/H \cong \mathfrak{m}$. Then*

$$R(c\partial_t + x, y; y, c\partial_t + x) = R^{\varphi(t)}(x, y; y, x) - \tfrac{1}{4}\Lambda^2 Q\left(\Lambda^2 \dot{\varphi}(x \wedge y), x \wedge y\right)$$

$$+ \; c \; \left(\tfrac{3}{2}Q\left(\dot{\varphi}[x,y], y\right) + Q\left(\varphi^{-1}\dot{\varphi}y, \pi^+(x,y)\right) - Q\left(\varphi^{-1}\dot{\varphi}x, \pi^+(y,y)\right)\right)$$

$$- \; \tfrac{1}{4}c^2 Q\left((2\ddot{\varphi} - \dot{\varphi}\varphi^{-1}\dot{\varphi})y, y\right),$$

*where $R^\varphi$ is the curvature of $(G/H, g_\varphi)$, $\pi^+(u,v) := \tfrac{1}{2}([u, \varphi(v)] + [v, \varphi(u)])$, $\Lambda^2 Q$ is the inner product on $\Lambda^2 \mathfrak{m}$ given as $\Lambda^2 Q(s \wedge t, u \wedge v) := Q(s,u)Q(t,v) - Q(s,v)Q(t,u)$, and where $\Lambda^2 \dot{\varphi} : \Lambda^2 \mathfrak{m} \to \Lambda^2 \mathfrak{m}$ is given by $\Lambda^2 \dot{\varphi}(u \wedge v) := (\dot{\varphi}u) \wedge (\dot{\varphi}v)$.*

We call a metric $g_\varphi$ *diagonalizable* and denote it by $g_f$ if there exists an $H$-invariant decomposition and a tuple of smooth positive functions

$$\mathfrak{m} = \mathfrak{m}_1 \oplus \ldots \oplus \mathfrak{m}_k, \quad f = (f_1, \ldots, f_k), \quad \text{such that} \quad \varphi(t)|_{\mathfrak{m}_i} = f_i(t)^2 Id_{\mathfrak{m}_i}. \tag{4}$$

We can now prove the following deformation result for positively curved metrics.

**Lemma 3.2** *Let $M := I \times G/H$ be endowed with a diagonalizable metric $g_f$ of the form (4) of positive sectional curvature. Moreover, suppose that at some $t_0 \in I$ we have $f_i(t_0) = c_0$ for all $i$.*

*Then there exists an $\varepsilon > 0$ such that on any open subinterval $t_0 \in J \subset I$ on which $|f_i - c_0| < \varepsilon$ and for any tuple of functions $\tilde{f} = (\tilde{f}_1, \ldots, \tilde{f}_k)$ on $J$ satisfying*

$$(i) \quad |\tilde{f}_i - c_0| < \varepsilon, \quad (ii) \; |\tilde{f}_i'| < |f_i'| + \varepsilon \quad (iii) \; \tilde{f}_i'' < f_i'' + \varepsilon$$

$$(iv) \quad \text{either } |\tilde{f}_i' - \tilde{f}_1'|, |f_i' - f_1'| < \varepsilon$$

$$\text{or } |\rho_k| < 1 + \varepsilon \quad \text{and} \quad |\rho_i - \rho_k| < \varepsilon \quad \text{where } \rho_i = \frac{\tilde{f}_i' - \tilde{f}_1'}{f_i' - f_1'}, \quad i = 2, \ldots, k$$

*the metric $(J \times G/H, g_{\tilde{f}})$ has also positive sectional curvature.*

**Proof.** We regard the curvature formula in Proposition 3.1 as a quadratic polynomial in $c$ and denote its coefficients by

$$R(c\partial_t + x, y; y, c\partial_t + x) = A(x,y) + B(x,y)\, c + C(y)\, c^2. \tag{5}$$

Note that $A$ and $C$ depend quadratically on $x \wedge y$ and $y$, respectively, whereas $B$ depends bilinearly on $x \wedge y$ and $y$. It follows that $(M, g_f)$ has positive sectional curvature iff for all $x, y \in \mathfrak{m}$ we have $A(x,y) > 0$ whenever $x \wedge y \neq 0$, $C(y) > 0$ whenever $y \neq 0$ and $A(x,y)C(y) - \tfrac{1}{4}B(x,y)^2 > 0$ whenever $x \wedge y \neq 0$.

Thus, the metric $(J \times G/H, g_{\tilde{f}})$ has also positive sectional curvature if the corresponding functions $\tilde{A}$, $\tilde{B}$ and $\tilde{C}$ satisfy $\tilde{A}(x,y) > A(x,y) - \delta ||x \wedge y||_Q^2$, $|B(x,y)| < |\tilde{B}(x,y)| + \delta ||x \wedge y||_Q ||y||_Q$ and $\tilde{C}(y) > C(y) - \delta ||y||_Q^2$ for some sufficiently small $\delta > 0$. Evidently, it suffices to show these inequalities in the case where $||y||_Q = ||x \wedge y||_Q = 1$, whence we shall assume this from now on.

Suppose now that $\delta > 0$ is given. Decomposing $x = \sum_i x_i$ and $y = \sum_i y_i$ according to (4), it follows from $\varphi|_{\mathfrak{m}_i} = f_i^2 Id_{\mathfrak{m}_i}$ that

$$\frac{1}{4}Q\left((2\ddot{\varphi} - \dot{\varphi}\varphi^{-1}\dot{\varphi})y, y\right) = \sum_i f_i f_i'' ||y_i||_Q^2.$$

Suppose for some $\varepsilon > 0$ we have functions $\tilde{f}_i$ as requested, and denote the corresponding metric by $\tilde{\varphi}$. Then on $J \subset I$ we have $|f_i - \tilde{f}_j| < 2\varepsilon$ for all $i, j$. Note that the eigenvalues of $\Lambda^2 \dot{\varphi}$ and $\Lambda^2 \dot{\tilde{\varphi}}$ are $4f_i f_j f_i' f_j'$ and $4\tilde{f}_i \tilde{f}_j \tilde{f}_i' \tilde{f}_j'$, respectively, and assuming that $\varepsilon > 0$ is sufficiently small, we may assume that $|\tilde{f}_i \tilde{f}_j \tilde{f}'_i \tilde{f}'_j| < |f_i f_j f'_i f'_j| + \delta/2$. Whence, for all $x, y \in \mathfrak{m}$ with $\|y\|_Q = \|x \wedge y\|_Q^2 = 1$ we get

$$\left| R^{\tilde{\varphi}(t)}(x, y; y, x) - R^{\varphi(t)}(x, y; y, x) \right| < \delta/2,$$

$$-\tfrac{1}{4}\Lambda^2 Q\left(\Lambda^2 \dot{\tilde{\varphi}}(x \wedge y), x \wedge y\right) > -\tfrac{1}{4}\Lambda^2 Q\left(\Lambda^2 \dot{\varphi}(x \wedge y), x \wedge y\right) - \delta/2 \quad \text{and}$$

$$-\tfrac{1}{4}Q\left((2\ddot{\tilde{\varphi}} - \dot{\tilde{\varphi}}\tilde{\varphi}^{-1}\dot{\tilde{\varphi}})y, y\right) = -\sum_i \tilde{f}_i \tilde{f}''_i \|y_i\|_Q^2 > -\sum_i (f_i f''_i + \delta)\|y_i\|_Q^2 = -\tfrac{1}{4}Q\left((2\ddot{\varphi} - \dot{\varphi}\varphi^{-1}\dot{\varphi})y, y\right) - \delta.$$

Therefore, $\tilde{A}(x, y) > A(x, y) - \delta$ and $\tilde{C}(y) > C(y) - \delta$.

Next, note that $B(x, y)$ depends polynomially on $f_i, 1/f_i$ and linearly on $f'_i$. Moreover, for $\varphi = \lambda Id$ and $\dot{\varphi} = \lambda' Id$ this expression vanishes; indeed, in this case $\pi^+(x, y) = \pi^+(y, y) = 0$ and $Q(\dot{\varphi}[x, y], y) = \lambda' Q([x, y], y) = \lambda' Q(x, [y, y]) = 0$ by the biinvariance of $Q$. Therefore, we must have

$$B(x, y) = \sum_{i>1} D_i(x, y, f_j) f'_1 (f_i - f_1) + \sum_{i>1} E_i(x, y, f_j)(f'_i - f'_1),$$

where $D_i$ and $E_i$ depend bilinearly on $x \wedge y$ and $y$ and polynomially on $f_j$, $1/f_j$.

Let us suppose that $|f'_i - f'_1|, |\tilde{f}'_i - \tilde{f}'_1| < \varepsilon$ for all $i \leq i_0$. Then, since $f_j, \tilde{f}_j$ are uniformly bounded on $J$, so are $D_i(x, y, f_j)$ and $E_i(x, y, f_j)$ for $\|y\|_Q = \|x \wedge y\|_Q = 1$. Thus, by choosing $\varepsilon > 0$ sufficiently small, we can achieve that

$$\left| B(x, y) - \sum_{i>i_0} E_i(x, y, f_j)(f'_i - f'_1) \right| < \frac{\delta}{4}, \quad \left| \tilde{B}(x, y) - \sum_{i>i_0} E_i(x, y, \tilde{f}_j)(\tilde{f}'_i - \tilde{f}'_1) \right| < \frac{\delta}{4}. \tag{6}$$

Now, we estimate

$$\left| \sum_{i>i_0} E_i(x, y, \tilde{f}_j)(\tilde{f}'_i - \tilde{f}'_1) \right| = \left| \sum_{i>i_0} E_i(x, y, \tilde{f}_j) \rho_i (f'_i - f'_1) \right| \leq \left| \sum_{i>i_0} E_i(x, y, f_j) \rho_i (f'_i - f'_1) \right| + \frac{\delta}{4}$$

$$\leq |\rho_k| \left| \sum_{i>i_0} E_i(x, y, f_j)(f'_i - f'_1) \right| + \sum_{i>i_0} |E_i(x, y, f_j)| |\rho_i - \rho_k| |f'_i - f'_1| + \frac{\delta}{4}$$

$$\leq \left| \sum_{i>i_0} E_i(x, y, f_j)(f'_i - f'_1) \right| + \frac{\delta}{2},$$

as long as $\varepsilon > 0$ is sufficiently small, and this together with (6) implies that $|\tilde{B}(x, y)| < |B(x, y)| + \delta$ which finishes the proof. ∎

**Proposition 3.3** *Let $M := I \times G/H$ be endowed with a diagonalizable metric $g_f$ of the form (4) of positive sectional curvature. Moreover, suppose that at some $t_0 \in I$ we have $f_i(t_0) = c_0$ and $f'_i(t_0) > 0$ for all $i$.*

*Then there exists a tuple of smooth functions $\tilde{f} = (\tilde{f}_1, \ldots, \tilde{f}_k)$ on $I = (a, b)$ such that $\tilde{f}_i \equiv f_i$ close to $a$ and $\tilde{f}_i \equiv f$ close to $b$, where $f : I \to \mathbb{R}$ is some smooth positive function with $f' > 0$, such that the metric $(I \times G/H, g_{\tilde{f}})$ has also positive sectional curvature.*

**Proof.** We let $\varepsilon > 0$ as in Lemma 3.2 and label the functions such that $f'_1(t_0) \leq \ldots \leq f'_k(t_0)$. By Lemma 3.2, we may replace $f_1$ by a $C^2$-close function and thus, after shrinking $\varepsilon$ if necessary, we may assume that $\varepsilon f'_k(t_0) < f'_1(t_0) < f'_2(t_0)$. By l'Hopital's rule, we have $\lim_{t \to t_0}(f_i(t) - f_1(t))/(f_k(t) - f_1(t)) = (f'_i(t_0) - f'_1(t_0))/(f'_k(t_0) - f'_1(t_0)) \in (0, 1]$ for $i > 1$, so that on some sufficiently small interval $t_0 \in J \subset I$ we have for all $s, t \in J$ the estimates

$$
\begin{aligned}
&(i) \quad |f_i(t) - c_0| < \frac{\varepsilon}{5} \quad \text{for } i \geq 1, \qquad (ii) \quad f_i'(t) > f_1'(t) > \varepsilon \; f_i'(t) > 0 \quad \text{for } i > 1, \\
&(iii) \quad 0 < \frac{f_i(t) - f_1(t)}{f_k(t) - f_1(t)} < 2 \quad \text{for } i > 1, \quad (iv) \quad N_i \left| \frac{f_i(t) - f_1(t)}{f_k(t) - f_1(t)} - \frac{f_i'(s) - f_1'(s)}{f_k'(s) - f_1'(s)} \right| < \varepsilon \quad \text{for } i > 1,
\end{aligned}
\tag{7}
$$

where
$$
N_i := \sup_{s \in J} \left| \frac{f_k'(s) - f_1'(s)}{f_i'(s) - f_1'(s)} \right| \quad \text{for } i > 1.
$$

Now it is easy to see that there exists a smooth function $\overline{f}_k : J = (a_1, b_1) \to \mathbb{R}$ satisfying

$$
\begin{aligned}
&(i) \; \overline{f}_k \equiv f_k \text{ close to } a_1, \qquad (ii) \; \overline{f}_k \leq f_1, f_k, \qquad (iii) \; |\overline{f}_k - c_0| < \frac{\varepsilon}{5}, \qquad (iv) \; \overline{f}_k^{(n)} \leq f_k^{(n)}, \; n = 1, 2, \\
&(v) \text{ for some } t_1 \in J, \; t_1 > t_0 \text{ we have } \overline{f}_k(t_1) = f_1(t_1), \; \overline{f}_k'(t_1) = f_1'(t_1) \text{ and } \overline{f}_k'(t) > f_1'(t) \text{ for all } t < t_1.
\end{aligned}
\tag{8}
$$

For $1 \leq i < k$, we define
$$
\overline{f}_i := f_i - C_i(f_k - \overline{f}_k), \qquad \text{where} \qquad C_i := \frac{f_i(t_1) - f_1(t_1)}{f_k(t_1) - f_1(t_1)}.
$$

Note that $C_1 = 0$ and hence $\overline{f}_1 = f_1$, and that $0 < C_i < 2$ for $i > 1$ by $(7)(iii)$. Moreover, from $(8)(i)$ and $(v)$ it follows that $\overline{f}_i \equiv f_i$ close to $a_1$, $\overline{f}_i(t_1) = f_1(t_1)$ and $\overline{f}_i'(t_1) = f_1'(t_1)$. If we let

$$
\rho_i := \frac{\overline{f}_i' - \overline{f}_1'}{f_i' - f_1'} = 1 - C_i \frac{f_k' - \overline{f}_k'}{f_i' - f_1'} \quad \text{for } i > 1,
$$

then using the estimates from (7) and (8) we have

$$
\begin{aligned}
&|\overline{f}_i(t) - c_0| \leq |f_i(t) - c_0| + C_i(|f_k(t) - c_0| + |c_0 - \overline{f}_k(t)|) < \varepsilon && \text{for all } i \geq 1 \text{ by } (7)(i) \; \& \; (8)(iii), \\
&\overline{f}_i^{(n)} = f_i^{(n)} - C_i(f_k^{(n)} - \overline{f}_k^{(n)}) \leq f_i^{(n)} < f_i^{(n)} + \varepsilon && \text{for } n = 1, 2, \text{ all } i \geq 1 \text{ by } (8)(iv), \\
&0 < \rho_k \leq 1 < 1 + \varepsilon && \text{for all } t \leq t_1 \text{ by } (8)(iv) \; \& \; (v), \\
&|\rho_i - \rho_k| = \left| \frac{f_i' - f_1'}{f_k' - f_1'} - C_i \right| \frac{|f_k' - \overline{f}_k'|}{|f_i' - f_1'|} \leq \left| \frac{f_i' - f_1'}{f_k' - f_1'} - C_i \right| N_i < \varepsilon && \text{for } t \leq t_1, \; i > 1 \text{ by } (8)(v) \; \& \; (7)(iv), \\
&\overline{f}_i' > \overline{f}_i' - f_1' + \varepsilon(f_i' - f_1') = (f_i' - f_1')(\rho_i + \varepsilon) && \text{by } (7) \; (ii) \\
&\qquad > (f_i' - f_1')\rho_k > 0 && \text{for } t \leq t_1 \text{ by the two preceding.}
\end{aligned}
\tag{9}
$$

In fact, after shrinking $J$ if necessary, we may assume that (9) holds on all of $J$. Now we pick a smooth function $f : I \to \mathbb{R}$ with $f' > 0$, $f'' < 0$, $f(t_1) = f_1(t_1)$, $f'(t_1) = f_1'(t_1)$ and $f''(t_1) < \overline{f}_i''(t_1)$ for all $i$. Then we can find smooth functions $\tilde{f}_i : I \to \mathbb{R}$ and some $\delta > 0$ such that

(i) $\tilde{f}_i \equiv f_i$ for $t \leq a_1$, \quad (ii) $\tilde{f}_i \equiv \overline{f}_i$ for $a_1 \leq t \leq t_1 - \delta$, \quad (iii) $\tilde{f}_i \equiv f$ for $t \geq t_1 + \delta$,

(iv) on $(t_1 - 2\delta, t_1 + 2\delta) \subset J$, we have $\tilde{f}_i'' \leq \overline{f}_i''$, and $|\tilde{f}_i - \overline{f}_i|, |\tilde{f}_i' - \overline{f}_i'|$ are sufficiently small such that (9) still holds when replacing $\overline{f}_i$ by $\tilde{f}_i$.

Then the metric $(M, g_{\tilde{f}})$ has positive sectional curvature; indeed, $g_{\tilde{f}} = g_f$ for $t \leq a_1$, and on $(a_1, t_1 + 2\delta)$, (9) holds so that we can use Lemma 3.2. Finally, for $t \geq t_1 + \delta$, $g_{\tilde{f}} = dt^2 + f^2 g_Q$. By Proposition 3.1, such a warped product has positive sectional curvature iff $f'' < 0$ and $(f')^2 < Sec(G/H, g_Q)$, thus the latter estimate holds for $t \in (t_1 + \delta, t_1 + 2\delta)$, and since $f' > 0$ and $f'' < 0$, it holds also for all $t \geq t_1 + 2\delta$. ∎

## 4 Applications of the Cheeger trick

In this section, we shall recall a well-known method which is very useful to construct metrics with lower sectional curvature bounds and goes back to J.Cheeger.

**Proposition 4.1** *[Ch] Let $(M, g)$ be a Riemannian manifold, let $(G, g_0)$ be a Lie group with Lie algebra $\mathfrak{g}$, equipped with a left invariant metric $g_0$, and suppose that $G$ acts isometrically on $M$ from the right. For $p \in M$, let $\mathfrak{h} \subset \mathfrak{g}$ be the infinitesimal stabilizer of $p$, and use the $g_0$-orthogonal decomposition $\mathfrak{g} = \mathfrak{h} \oplus \mathfrak{m}$. Let*

$$T_p M = V_1 \oplus V_2 \cong \mathfrak{m} \oplus V_2$$

*be the $g$-orthogonal decomposition, where $V_1 := T_p(G \cdot p) \cong \mathfrak{m}$, denoting the latter identification by $x_* \in T_p(G \cdot p) \leftrightarrow x \in \mathfrak{m}$. Let $\varphi_p : \mathfrak{m} \to \mathfrak{m}$ be the endomorphism defined by*

$$g_p(x_*, y_*) = g_0(\varphi_p x, y) \quad \text{for all } x, y \in \mathfrak{m}.$$

*Then there is a Riemannian metric $g'$ on $M$ for which the multiplication map $(M \times G, g + g_0) \to (M, g')$ is a Riemannian submersion, and $g'$ is determined by*

$$g'(x_*, y_*) = g_0(\varphi_p(\varphi_p + Id)^{-1}x, y) \quad \text{for all } x, y \in \mathfrak{m}, \quad g'|_{V_2} = g|_{V_2}, \quad g'(V_1, V_2) = 0.$$

*Furthermore, suppose that $(G, g_0)$ has nonnegative sectional curvature. Then $(M, g')$ has nonnegative (positive) sectional curvature provided $(M, g)$ has nonnegative (positive) sectional curvature.*

**Proof.** Fix $p \in M$, and note that the fiber of the map $\pi : M \times G \to M$ consists of all pairs $(p \cdot g, g^{-1})$, whence the vertical and horizontal subspaces are

$$\mathcal{V}_p := (0, \mathfrak{h}) \oplus \{(x_*, -x) \mid x \in \mathfrak{m}\}, \quad \mathcal{H}_p := (V_2, 0) \oplus \{(x_*, \varphi_p x) \mid x \in \mathfrak{m}\}.$$

Indeed, $(g_0 + g)(\mathcal{H}_p, \mathcal{V}_p) = 0$ is easily verified. Note that $\pi_*(x_*, \varphi_p x) = (\varphi_p + Id)x_*$, whence the horizontal lift of $x_* \in V_1$ is $((\varphi_p + Id)^{-1}x_*, \varphi_p(\varphi_p + Id)^{-1}x)$, whereas the horizontal lift of $v \in V_2$ is $(v, 0)$. Thus,

$$\begin{aligned}
g'(x_*, y_*) &= g((\varphi_p + Id)^{-1}x_*, (\varphi_p + Id)^{-1}y_*) + g_0(\varphi_p(\varphi_p + Id)^{-1}x, \varphi_p(\varphi_p + Id)^{-1}y) \\
&= g_0(\varphi_p(\varphi_p + Id)^{-2}x, y) + g_0(\varphi_p^2(\varphi_p + Id)^{-2}x, y) \\
&= g_0(\varphi_p(\varphi_p + Id)^{-1}x, y).
\end{aligned}$$

The remaining equations are evident, and the final statement is an immediate consequence of O'Neill's formula and the observation that the horizontal lift of a tangent plane of $M$ is always transversal to $\{p\} \times G \subset M \times G$. ∎

We shall use this method to obtain the following result on *cohomogeneity one manifolds*, i.e. manifolds $M$ on which a compact Lie group $G$ acts smoothly with principal orbits of codimension one (cf. Section 6).

**Theorem 4.2** *[ST] Let $(M, G)$ be a (closed) cohomogeneity one manifold. Then $M$ has $G$-invariant metrics of almost nonnegative sectional curvature, i.e. for every $\varepsilon > 0$ there is an invariant metric $g_\varepsilon$ on $M$ such that $Sec(M, g_\varepsilon) \cdot diam(M, g_\varepsilon)^2 > -\varepsilon$.*

The original proof in [ST] gave explicit constructions of the metrics $g_\varepsilon$; indeed, this result also follows from Theorem 6.1. However, the observation that the Cheeger trick allows to obtain this result much faster is due to B.Wilking, and we shall present his proof here.

**Proof.** It is well known that any cohomogeneity one manifold which has less than two singular orbits admits invariant metrics of nonnegative sectional curvature, hence we can restrict to the case where $M$ is closed and has two singular orbits.

We pick a $G$-invariant metric $g$ on $M$ which has nonnegative sectional curvature in a neighborhood of the singular orbits. Moreover, let $C > 0$ be a constant such that

$$Sec(M, g) \geq -C.$$

Let $Q$ be a biinvariant metric on $G$, and define for $\delta > 0$ the metric $g_\delta$ on $M$ by the property that the submersion $(G \times M, \delta Q + g) \to (M, g_\delta)$ is Riemannian. Since Riemannian submersions do not increase the diameter, we have for $\delta \leq 1$

$$diam(M, g_\delta) \leq diam(G \times M, \delta Q + g) \leq diam(G \times M, Q + g) =: D,$$

so that this diameter is uniformly bounded. For $p \in M$ and $\delta \in (0, 1]$ we define the endomorphisms $\varphi_p^\delta : \mathfrak{m} \to \mathfrak{m}$ by the equation

$$g_p(x_*, y_*) = \delta Q\left(\varphi_p^\delta x, y\right) \quad \text{for all } x, y \in \mathfrak{m},$$

so that evidently $\varphi_p^\delta = \delta^{-1} \varphi_p^1 =: \delta^{-1} \varphi_p$. Thus, if $p \in M$ is contained in a principal orbit then according to Proposition 4.1 the horizontal lift of $c\partial_t + x_* \in T_p M$ is

$$\overline{c\partial_t + x_*} = \left(\varphi_p^\delta(\varphi_p^\delta + Id)^{-1}x, c\partial_t + (\varphi_p^\delta + Id)^{-1}x_*\right)$$

$$= \left(\varphi_p(\varphi_p + \delta Id)^{-1}x, c\partial_t + \delta(\varphi_p + \delta Id)^{-1}x_*\right) \in T_{(e,p)}(G \times M).$$

Now a straightforward calculation yields

$$\lim_{\delta \to 0} \delta^{-1} \|(c\partial_t + x_*) \wedge y_*\|_{g_\delta}^2 = c^2 \|y\|_Q,$$

$$\lim_{\delta \to 0} \delta^{-2} \|x_* \wedge y_*\|_{g_\delta}^2 = \|x \wedge y\|_Q^2.$$

On the other hand, by O'Neill's formula and the fact that $(G, \delta Q)$ has nonnegative sectional curvature, we have

$$R_M^{g_\delta}((c\partial_t + x_*) \wedge y_*) \geq R_M^g((c\partial_t + \delta(\varphi_p + \delta Id)^{-1}x_*) \wedge \delta(\varphi_p + \delta Id)^{-1}y_*)$$

$$\geq -C\delta^2 \|(c\partial_t + \delta(\varphi_p + \delta Id)^{-1}x_*) \wedge (\varphi_p + \delta Id)^{-1}y_*\|_g^2,$$

so that

$$\liminf_{\delta \to 0} \delta^{-2} R_M^{g_\delta}((c\partial_t + x_*) \wedge y_*) \geq -Cc^2 \|(\varphi_p)^{-1}y_*\|_g^2,$$

$$\liminf_{\delta \to 0} \delta^{-4} R_M^{g_\delta}(x_* \wedge y_*) \geq -C\|(\varphi_p)^{-1}x_* \wedge (\varphi_p)^{-1}y_*\|_g^2.$$

Taking the quotient, it follows that in either case

$$\liminf_{\delta \to 0} \delta^{-1} Sec_{g_\delta}((c\partial_t + x_*) \wedge y_*) > -\infty, \quad \text{whence} \quad \liminf_{\delta \to 0} Sec_{g_\delta}((c\partial_t + x_*) \wedge y_*) \geq 0.$$

By compactness of $M$, we may assume that this convergence is uniform, and since $diam(M, g_\delta) \leq D$, we have that $\liminf_{\delta \to 0} Inf(Sec(M, g_\delta)) \cdot diam(M, g_\delta)^2 \geq 0$ which proves the claim. ∎

We shall use another well known application of the Cheeger trick. Namely, suppose that we have a chain of compact Lie groups $K_0 \subset \ldots \subset K_n = G$, and fix a biinvariant metric $Q$ on $\mathfrak{g}$. We decompose the Lie algebra $\mathfrak{g}$ $Q$-orthogonally as

$$\mathfrak{g} = \mathfrak{m}_0 \oplus \ldots \oplus \mathfrak{m}_n, \quad \mathfrak{k}_i = \mathfrak{m}_0 \oplus \ldots \oplus \mathfrak{m}_i. \tag{10}$$

For constants $c_0, \ldots, c_n > 0$ we define an inner product $g = g_{(c_0,\ldots,c_n)}$ on $\mathfrak{g}$ by

$$g(\mathfrak{m}_i, \mathfrak{m}_j) = 0 \quad \text{for all } i \neq j, \quad g|_{\mathfrak{m}_i} = c_i Q|_{\mathfrak{m}_i} \quad \text{for all } i, \tag{11}$$

and denote the corresponding left invariant metric on $G$ by the same symbol. Then we have

**Proposition 4.3** *Let $G$ be a compact Lie group as above and define the left invariant metric $g$ as in (11) using the decomposition (10). If $0 < c_0 \leq \ldots \leq c_n$ then $g_{(c_0,\ldots,c_n)}$ has nonnegative sectional curvature.*

**Proof.** We proceed by induction on $n$. For $n = 0$, $g = c_0 Q$ which has nonnegative sectional curvature.

In general, given $K_0 \subset \ldots \subset K_n = G$, we know by induction hypothesis that for all constants $c_1 \leq \ldots \leq c_n$ the metric $g_{(c_1,c_1,\ldots,c_n)}$ has nonnegative sectional curvature, since it is associated to the chain $K_1 \subset \ldots \subset K_n = G$. This already establishes the claim in the case $c_0 = c_1$.

If $c_0 < c_1$ then we apply the Cheeger trick to the right action of $K_0$ on $(G, g_{(c_1,c_1,\ldots,c_n)})$, where we equip $K_0$ with the (nonnegatively curved) metric $\lambda Q$, some $\lambda > 0$. Then by Proposition 4.1, the resulting metric $g'$ has nonnegative sectional curvature, and by the formula given there, $g' = g_{(c_0, c_1,\ldots,c_n)}$, where $c_0 = \frac{\lambda}{1+\lambda} c_1$. Evidently, for the appropriate choice of $\lambda > 0$, we obtain any $c_0 \in (0, c_1)$. ∎

With this, we now obtain the following result.

**Proposition 4.4** *Let $K \subset O(n+1)$ be a Lie subgroup which acts transitively on $S^n \subset \mathbb{R}^{n+1}$. Let $Q$ be a biinvariant metric on $\mathfrak{k}$, and denote the induced normal homogeneous metric on $S^n = K/H$ by $g_Q$.*

*Then there is a chain $H = K_0 \subset \ldots \subset K_r = K$ and a left invariant metric $g_K = g_{(c_0,\ldots,c_r)}$ of nonnegative sectional curvature on $K$ such that the submersion $(K \times S^n, g_K + g_0) \to (S^n, \mu g_Q)$ is Riemannian for some constant $\mu > 0$, where $g_0$ denotes the standard metric on $S^n$.*

**Proof.** Transitive actions on spheres have been extensively studied in the literature. From their classification ([MS]) it follows that there is always a chain of Lie groups $H = K_0 \subset \ldots \subset K_r = K$ such that when using the decomposition $\mathfrak{g} = \mathfrak{h} \oplus \mathfrak{m}_1 \oplus \ldots \oplus \mathfrak{m}_r$ from (10), there are constants $\rho_1, \ldots, \rho_r > 0$ such that the submersion $(K, g_{(\rho_0,\rho_1,\ldots,\rho_r)}) \to (S^n, g_0)$ is Riemannian. Note that $\rho_1 \geq \ldots \geq \rho_r > 0$ (cf. Table 2.4 in [GZ2]) and that $\rho_0 > 0$ can be chosen arbitrarily, so that we may assume $\rho_0 \geq \rho_1$.

By Proposition 4.1 it follows that if we fix $\mu \in (0, \rho_r)$ and let $c_i := \mu \rho_i / (\rho_i - \mu)$, then $0 < c_0 \leq \ldots \leq c_r$, and there are Riemannian submersions

$$(K \times K, g_{(c_0,\ldots,c_r)} + g_{(\rho_0,\ldots,\rho_r)}) \to (K, \mu Q), \quad \text{which induces} \quad (K \times S^n, g_{(c_0,\ldots,c_r)} + g_0) \to (S^n, \mu g_Q),$$

and $(K, g_{(c_0,\ldots,c_r)})$ has nonnegative sectional curvature by Proposition 4.3. ∎

**Theorem 4.5** *Let $K \subset O(n+1)$ be a Lie subgroup which acts transitively on $S^n \subset \mathbb{R}^{n+1}$, and let $g_Q$ be a normal homogeneous metric on $S^n$ induced by some $\mathrm{Ad}_K$-invariant inner product $Q$ on $\mathfrak{k}$. Let $r(x) := \|x\|$ be the radius function on $\mathbb{R}^{n+1}$.*

*Then there exists a $K$-invariant metric $g$ on the unit ball $B_1(0) \subset \mathbb{R}^{n+1}$ with positive sectional curvature, and an $\varepsilon > 0$, such that on $r^{-1}(1-\varepsilon, 1)$ we have $g = dr^2 + f(r)^2 g_Q$ where $f : (1-\varepsilon, 1) \to \mathbb{R}$ satisfies $f > 0, f' > 0$.*

**Proof.** Let us consider a (complete) $O(n+1)$-invariant metric $\tilde{g}$ of positive sectional curvature on $\mathbb{R}^{n+1}$. By the $O(n+1)$-invariance, the restriction of $\tilde{g}$ to each sphere centered at $0$ must be a multiple of the round metric $g_0$, hence in polar coordinates on $\mathbb{R}^{n+1} \setminus 0 \cong (0, \infty) \times K/H$ we may write

$$\tilde{g} = dt^2 + \lambda(t)^2 g_0$$

for some function $\lambda(t) : (0, \infty) \to \mathbb{R}$. We may choose $\tilde{g}$ such that for some $t_0 \in (0,1)$ we have $\lambda(t_0) = 1$ and $\lambda'(t_0) > 0$. Let $g_K$ be the left invariant metric of nonnegative sectional curvature on $K$ from Proposition 4.4 and consider the $K$-invariant metric $\overline{g}$ on $\mathbb{R}^{n+1}$ for which the submersion

$$(K \times \mathbb{R}^{n+1}, g_K + \tilde{g}) \longrightarrow (\mathbb{R}^{n+1}, \overline{g})$$

is Riemannian. By Proposition 4.1 it follows that $\overline{g}$ has also positive sectional curvature and is diagonal of the form

$$\overline{g} = dt^2 + g_{(f_1,\ldots,f_r)}, \quad \text{where} \quad f_i(t)^2 = \frac{\lambda(t)\mu\rho_i}{\lambda(t)\rho_i + (1-\lambda(t))\mu}.$$

Evidently, $f_i(t_0) = \sqrt{\mu}$ and $f_i'(t_0) > 0$ for all $i$. Thus, by Proposition 3.3 we can construct a $K$-invariant metric $g$ of positive sectional curvature on $B_1(0) \backslash 0 = (0,1) \times S^n$ which close to the origin coincides with $\overline{g}$ and hence can be extended to 0, and which close to the boundary of $B_1(0)$ is of the asserted form. ∎

## 5  Homogeneous vector bundles

Let $G/K$ be a compact homogeneous space, and suppose there is a representation $\iota : K \to Aut(V)$ on some finite dimensional vector space $V$, which we may assume to be orthogonal as $K$ is compact. Then we can associate the *homogeneous vector bundle*

$$D := G \times_K V,$$

i.e. the set of equivalence classes under the relation on $G \times V$ given by $(gh, v) \sim (g, hv)$ for all $g \in G$, $h \in K$ and $v \in V$. Thus, we can regard $D$ as the orbit space of $G \times V$ under the "diagonal action" $h \cdot (g, v) := (gh^{-1}, hv)$ of $K$, and since this action is free, it follows that for any $K$-invariant metric on $G \times V$ we get a (unique) metric on $D$ for which the submersion $G \times V \to D$ is Riemannian.

Note that there is a canonical action of $G$ on $D$, and the cohomogeneity of the principal orbit of this action equals the cohomogeneity of the principal orbit of the action of $K$ on $V$.

Let us assume that $G$ and hence $K$ act by cohomogeneity one. Since $K$ acts orthogonally, it leaves all spheres centered at the origin invariant, whence this action has cohomogeneity one iff $K$ acts transitively on the unit sphere $S^n \subset V \cong \mathbb{R}^{n+1}$. In particular, we can write the unit sphere $S^n = K/H$ as a homogeneous space where $H \subset K$ is the stabilizer of some unit vector in $V$.

Note that the norm function $r : V \to \mathbb{R}$, $v \mapsto ||v||$ is $K$-invariant and hence induces a function $r : D \to \mathbb{R}$, and for $R \in \mathbb{R}$, we let

$$D_R := r^{-1}([0, R)) \subset D. \tag{12}$$

Moreover, the level sets of $r$ are precisely the $G$-orbits of $D$. We fix a biinvariant inner product $Q$ on $\mathfrak{g}$, and choose subspaces $\mathfrak{m}_1, \mathfrak{m}_2 \subset \mathfrak{g}$ such that

$$\mathfrak{g} = \mathfrak{h} \oplus \mathfrak{m}_1 \oplus \mathfrak{m}_2 \quad \text{and} \quad \mathfrak{k} = \mathfrak{h} \oplus \mathfrak{m}_1 \tag{13}$$

are $Q$-orthogonal decompositions. We consider the metric on $B_1(0) \subset V \cong \mathbb{R}^{n+1}$ from Theorem 4.5. Evidently, we can extend this to a nonnegatively curved $K$-invariant metric $g_V$ on all of $V$ which outside of $B_1(0)$ has the form

$$g_V = dt^2 + c_0^2 t^2 g_Q \quad \text{for some constant } c_0 > 0.$$

The corresponding submersion metric on $D$ has nonnegative sectional curvature, and by Proposition 4.1, it can be written on $D \backslash D_1$ in the form

$$g = dt^2 + f_0(t)^2 Q|_{\mathfrak{m}_1} + Q|_{\mathfrak{m}_2}, \quad \text{where} \quad f_0(t) := \frac{c_0 t}{\sqrt{1 + c_0^2 t^2}}. \tag{14}$$

In general, if a metric on $I \times G/H$ is of the form

$$g = dt^2 + f(t)^2 Q|_{\mathfrak{m}_1} + Q|_{\mathfrak{m}_2}, \tag{15}$$

for some smooth function $f > 0$, then the curvature of $g$ can be calculated from Proposition 3.1 as follows.

**Proposition 5.1** *[ST] [GZ2] Consider the metric $g$ on $I \times G/H$ as in (15), and let $c \in \mathbb{R}$ and $x = x_1 + x_2, y = y_1 + y_2 \in \mathfrak{m}$ with $x_i, y_i \in \mathfrak{m}_i$. Using the identification $\mathfrak{m} \cong T_{eH} G/H$ and denoting the unit tangent vector of $I$ by $\partial_t$, the curvature of $(I \times G/H, g)$ is given by*

$$R(c\partial_t + x, y; y, c\partial_t + x) =$$
$$\tfrac{3}{4} f^2 \|[x,y]_\mathfrak{h}\|_Q^2 + \tfrac{1}{4} \|[x_2, y_2]_2 + f^2 ([x_1, y_2] + [x_2, y_1])\|_Q^2$$
$$+ \tfrac{1}{4} f^2 \|[x_1, y_1]\|_Q^2 + \tfrac{1}{2} f^2 (3 - 2f^2) Q([x_1, y_1], [x_2, y_2]) + (1 - \tfrac{3}{4} f^2) \|[x_2, y_2]_\mathfrak{k}\|_Q^2 \quad (16)$$
$$+ 3cff' Q([x_2, y_2], y_1) - c^2 ff'' Q(y_1, y_1) - (ff')^2 \|x_1 \wedge y_1\|_Q^2,$$

*where the subscripts refer to the decompositions (13).*

From there, one can now deduce the following

**Corollary 5.2** *Consider the metric $g$ on $I \times G/H$ as in (15), and let $c \in \mathbb{R}$ and $x = x_1 + x_2, y = y_1 + y_2 \in \mathfrak{m}$ with $x_i, y_i \in \mathfrak{m}_i$. If $0 < f \leq 1$ satisfies the inequality*

$$-ff'' \geq 9(f')^2, \quad (17)$$

*then*

$$R(c\partial_t + x, y; y, c\partial_t + x) \geq -(f'f)^2 \|x_1 \wedge y_1\|_Q^2 = -\left(\frac{f'}{f}\right)^2 \|x_1 \wedge y_1\|_g^2.$$

*In particular, $\mathrm{Sec}(I \times G/H, g) \geq -\left(\frac{f'}{f}\right)^2$. Moreover, if $f' \neq 0$ and $x_1 = 0$, then*

$$R(c\partial_t + x_2, y; y, c\partial_t + x_2) = 0 \text{ iff } cy_1 = 0 \text{ and } [x_2, y_2] + f^2 [x_2, y_1] = 0.$$

**Proof.** We decompose $[x_2, y_2]_\mathfrak{k} = v + w$ where $Q(v, y_1) = 0$ and $w$ is a scalar multiple of $y_1$. Since $Q([x_1, y_1], y_1) = Q(x_1, [y_1, y_1]) = 0$, it follows that $Q([x_1, y_1], [x_2, y_2]) = Q([x_1, y_1], v)$ and $Q(y_1, [x_2, y_2]) = Q(y_1, w)$. Finally, $\|[x_2, y_2]_\mathfrak{k}\|_Q^2 = \|v\|_Q^2 + \|w\|_Q^2$. Substituting all of this, (16) becomes

$$R(c\partial_t + x, y; y, c\partial_t + x) =$$
$$\tfrac{3}{4} f^2 \|[x,y]_\mathfrak{h}\|_Q^2 + \tfrac{1}{4} \|[x_2, y_2]_2 + f^2 ([x_1, y_2] + [x_2, y_1])\|_Q^2$$
$$+ \quad \tfrac{1}{4} f^2 Q([x_1, y_1], [x_1, y_1]) + \tfrac{1}{2} f^2 (3 - 2f^2) Q([x_1, y_1], v) + (1 - \tfrac{3}{4} f^2) Q(v, v)$$
$$+ \quad (1 - \tfrac{3}{4} f^2) Q(w, w) + 3ff' Q(w, cy_1) - ff'' Q(cy_1, cy_1)$$
$$- \quad (ff')^2 \|x_1 \wedge y_1\|_Q^2.$$

The first line of the right hand side of this equation is evidently nonnegative; the second and third lines are nonnegative since the quadratic polynomials

$$\tfrac{1}{4} f^2 x^2 + \tfrac{1}{2} f^2 (3 - 2f^2) x + (1 - \tfrac{3}{4} f^2) \quad \text{and} \quad (1 - \tfrac{3}{4} f^2) x^2 + 3ff' x - ff''$$

are nonnegative for all $x$; indeed, their leading coefficients are positive, and their discriminants are

$$\tfrac{1}{4} f^2 (1 - \tfrac{3}{4} f^2) - \tfrac{1}{16} f^4 (3 - 2f^2)^2 = \tfrac{1}{4} f^2 (1 - f^2)^3 \geq 0 \quad \text{and}$$

$$-ff''(1 - \tfrac{3}{4} f^2) - \tfrac{9}{4} (ff')^2 \geq 9(f')^2 (1 - \tfrac{3}{4} f^2) - \tfrac{9}{4} (ff')^2 = 9(f')^2 (1 - f^2) \geq 0.$$

Moreover, if $f' \neq 0$ then $f < 1$ as 1 is the maximal value for $f$, so that these discriminants are positive. Thus, equality holds only under the conditions state above. ∎

Let $L \subset G$ be a closed Lie subgroup with Lie algebra $\mathfrak{l} \subset \mathfrak{g}$ which acts freely on $G/K$ and hence acts also freely on $D$ and $G/H$. We denote the quotient by $\underline{D} := L\backslash D \to \Sigma$ where $\Sigma := L\backslash G/K$, and this is a vector bundle of the same rank as $D \to G/K$. Note that

$$D\backslash 0 = G \times_K (\mathbb{R}^+ \times K/H) \cong \mathbb{R}^+ \times G/H, \quad \text{and} \quad \underline{D}\backslash 0 \cong \mathbb{R}^+ \times L\backslash G/H.$$

Let us now suppose that for some $R_0$ the metric on $D\backslash D_{R_0} = (R_0, \infty) \times G/H$ is of the form (15) using the decomposition (13). Then the vertical and horizontal subspace of the submersion $\pi : D \to \underline{D}$ at a point $(t, gH)$ is given by

$$\mathcal{V}_g := pr_{\mathfrak{m}}(Ad_{g^{-1}}\mathfrak{l}), \quad \text{and} \quad \mathcal{H}_g := \mathbb{R}\partial_t \oplus \{y_1 + f(t)^2 y_2 \mid y_i \in \mathfrak{m}_i,\ y_1 + y_2 \in \mathcal{H}_g^Q\},$$

with

$$\mathcal{H}_g^Q := (\mathcal{V}_g)_Q^\perp,$$

where $pr_{\mathfrak{m}} : \mathfrak{g} \to \mathfrak{m}$ is the $Q$-orthogonal projection, and where we use again the identification $T_{gH}G/H \cong \mathfrak{m}$.

**Theorem 5.3** *Let $H \subset K \subset G$ and $Q$ be as before, and consider a metric $g$ on $M := I \times G/H$ as in (15) where $f : I \to (0, 1]$ is a smooth function. Let $L \subset G$ be a group which acts freely on $G/K$ and hence on $G/H$ and $M$, and denote the metric on $N := I \times L\backslash G/H$ induced by the canonical submersion $\pi : M \to N$ by $\underline{g}$. Evidently, $\underline{g}$ is invariant under the action of $Norm(L)/L$ on $N$. If $f$ satisfies*

$$-ff'' \geq C(f')^2, \quad \text{or, equivalently,} \quad (f^{C+1})'' \leq 0, \quad \text{where} \quad C := \max\{\dim N, 9\}, \tag{18}$$

*then $(N, \underline{g})$ has nonnegative Ricci curvature.*

*Moreover, if $\dim K > \dim H$ and if the biquotient $\Sigma := L\backslash G/K$ has finite fundamental group, then for all $t_0 \in I$ with $f'(t_0) \neq 0$, there is a $p \in L\backslash G/H$ with $Ric(N, \underline{g}) > 0$ at $(t_0, p) \in N$.*

**Proof.** By (16), we have $R_M(\partial_t, x; x, \partial_t) = -ff''\|x_1\|_Q^2$, and we choose a $g$-orthonormal basis $\partial_t, y^1, \ldots, y^k$ of $\mathcal{H}_g$ where $k = \dim N - 1$. Since (18) implies (17), Corollary 5.2 yields

$$R_M(c\partial_t + x, y^i; y^i, c\partial_t + x) \geq -(f'f)^2\|x_1 \wedge y_1^i\|_Q^2 \geq -(f'f)^2\|x_1\|_Q^2\|y_1^i\|_Q^2 = -(f')^2\|x_1\|_Q^2,$$

as $f^2\|y_1^i\|_Q^2 = \|y_1^i\|_g^2 = 1$. Thus, O'Neill's formula and (18) imply

$$\begin{aligned} Ric_N(c\partial_t + \underline{x}) &= R_N(\partial_t, \underline{x}; \underline{x}, \partial_t) + \sum_{i=1}^k R_N(c\partial_t + \underline{x}, y^i; y^i, c\partial_t + \underline{x}) \\ &\geq R_M(\partial_t, x; x, \partial_t) + \sum_{i=1}^k R_M(c\partial_t + x, y^i; y^i, c\partial_t + x) \\ &\geq \left(-ff'' - k(f')^2\right)\|x_1\|_Q^2 \geq (C - \dim N + 1)(f')^2\|x_1\|_Q^2 \geq 0. \end{aligned} \tag{19}$$

If $Ric_N(c\partial_t + \underline{x}) = 0$ and $f'(t_0) \neq 0$, then (19) implies that $x_1 = 0$, i.e. $x \in \mathfrak{m}_2 \cap \mathcal{H}_g = \mathfrak{m}_2 \cap \mathcal{H}_g^Q$, in which case Corollary 5.2 implies that $R_M(\partial_t, x; x, \partial_t) = 0$ and $R_M(c\partial_t + x, y^i; y^i, c\partial_t + x) \geq 0$. Thus, we must have that $R_M(c\partial_t + x, y^i; y^i, c\partial_t + x) = 0$ for all $i$ which, again by Corollary 5.2, implies that $cy_1^i = 0$ and $[x, y_1^i + f^2 y_2^i] = 0$ for all $i$. Since $\mathfrak{m}_1 \neq 0$ and $L$ acts freely on $G/K$ so that $\mathcal{V}_g \cap \mathfrak{m}_1 = 0$, there must be some $i$ for which $y_1^i \neq 0$, hence the first condition implies that $c = 0$.

On the other hand, $y_1^i + y_2^i \in \mathcal{H}_g$ iff $y_1^i + f^2 y_2^i \in \mathcal{H}_g^Q$ so that the second condition implies that $[x, \mathcal{H}_g^Q] = 0$. Note, however, that $\mathfrak{m}_2 \cap \mathcal{H}_g^Q$ is the tangent space of the biquotient $\Sigma$, so that this condition implies that $Ric_\Sigma^{g_Q}(x) = 0$. One the other hand, by Theorem 2.1, $Ric_\Sigma^{g_Q} > 0$ on an open dense subset of $\Sigma$, and this completes the proof. ∎

**Corollary 5.4** *Let $D := G \times_K V \to G/K$ be a homogeneous vector bundle where $K \subset O(V)$ acts transitively on the unit sphere $S^n = K/H$, so that we have the inclusions $H \subset K \subset G$. Moreover, suppose that $L \subset G$ acts freely on $G/K$ and hence on $G/H$ and $D$. For $R > 0$, we define $D_R$ as in (12).*

*Then for every $\delta > 0$ there is a metric $g_\delta$ on $\underline{D}_R := L \backslash D_R$ with the following properties:*

1. *$Sec(\underline{D}_R, g_\delta) \geq -\delta$ and $diam(\underline{D}_R, g_\delta) = O(\delta^{-1/6})$.*

2. *$Ric(\underline{D}_R, g_\delta) \geq 0$.*

3. *If $\dim V \geq 2$ and $L \backslash G/K$ has finite fundamental group, then $(\underline{D}_R, g_\delta)$ contains points of positive Ricci curvature.*

4. *The action of $Norm_G(L)/L$ on $\underline{D}_R$ is isometric.*

5. *Close to the boundary, $g_\delta$ is isometric to $(I \times L \backslash G/H, dt^2 + g_Q)$ where $g_Q$ stands for the metric on $L \backslash G/H$ induced by the biinvariant metric on $G$.*

Here, $O(\delta^p)$ denotes a function such that $\limsup_{\delta \to 0} |\delta^{-p} O(\delta^p)| < \infty$.

**Proof.** Recall the metric on $D$ from (14) which outside of $D_1$ has the form

$$g = dt^2 + f_0^2 Q|_{\mathfrak{m}_1} + Q|_{\mathfrak{m}_2}, \quad \text{where} \quad f_0(t) := \frac{c_0 t}{\sqrt{1 + c_0^2 t^2}}.$$

We define $R_0 = R_0(\delta)$ implicitly by

$$\left(\frac{f_0'}{f_0}\right)^2 (R_0) = \delta,$$

and verify that $R_0 = O(\delta^{-1/6})$. Moreover, if $\delta$ is sufficiently small so that $R_0$ is large, one verifies that $\left(f_0^{C+1}\right)''(R_0) < 0$ where $C := \max\{\dim \underline{D}_R, 9\}$ is defined as in (18). Thus, if we let

$$R := R_0 + \frac{1 - f_0^{C+1}}{(f_0^{C+1})'}(R_0) + 1$$

then it is straightforward to verify that $R = O(\delta^{-1/6})$, and there is a smooth function $h : [1, R] \to (0, 1]$ such that

$$h|_{[1, R_0]} = f_0^{C+1}, \quad h''|_{[R_0, R)} \leq 0, \quad h \equiv 1 \text{ close to } R. \tag{20}$$

Now we let $f := h^{1/(C+1)}$ and define the metric $g_\delta$ on $D$ by extending the metric on $D_1$ from above by the metric $g_f$ on $D \backslash D_1 = (1, \infty) \times G/H$. Since $f = f_0$ on $[1, R_0]$, it follows that on $D_{R_0}$, $g_\delta$ coincides with $g$ and hence has nonnegative sectional curvature. Also, (20) implies that $f'' \leq 0$ on $[R_0, R)$ and thus $f'/f$ decreases. As $f'/f(R_0) = \sqrt{\delta}$, we have $Sec(D_R \backslash D_{R_0}, g_\delta) \geq -\delta$ by Corollary 5.2. Combining these two estimates, we have $Sec(D_R, g_\delta) \geq -\delta$ and hence the same estimate holds for the submersion metric on $\underline{D}_R$.

Also, $g_\delta$ has nonnegative sectional and thus nonnegative Ricci curvature on $\underline{D}_{R_0}$. Moreover, by (20) and Theorem 5.3, $Ric(\underline{D}_R \backslash \underline{D}_{R_0}, g_\delta) \geq 0$ as well, and there exist points of positive Ricci curvature on $\underline{D}_R$ under the asserted conditions.

Further, close to the boundary of $D_R$, $g_\delta$ is the product $dt^2 + g_Q$ where $g_Q$ stands for the metric induced by the biinvariant one. Thus, the same is true for the submersion metric close to the boundary of $\underline{D}_R$.

Finally, $g_\delta|_{D_1}$ is independent of $\delta$ and since on $D \backslash D_1$, $g_\delta$ is of the form $dt^2 + g_f$, we have $diam(D_R, g_\delta) \leq 2R + diam(D_1) = O(\delta^{-1/6})$, and since Riemannian submersions do not increase the diameter, this completes the proof. ∎

# 6 Quotients of cohomogeneity one manifolds

Let $(M, G)$ be a cohomogeneity one manifold. It is well known that $M$ has at most two singular orbits. Indeed, if there are no singular orbits, then $M$ is a homogeneous bundle over either $\mathbb{R}$ or $S^1$, and in either case one easily constructs $G$-invariant metrics of nonnegative sectional curvature on $M$. If $M$ has one singular orbit, then $M$ is a homogeneous disc bundle of the type considered in the preceding section and hence also admits $G$-invariant metrics of nonnegative sectional curvature.

Finally, if $M$ has two singular orbits, then $M$ can be obtained by glueing together two disc bundles of the form $D_R$ considered above. Thus, if $L \subset G$ acts freely on $M$, then the quotient $L \backslash M$ is obtained by glueing two disc bundles of the form $\underline{D}_R$ along their common boundary. Therefore, we obtain as an immediate consequence of Corollary 5.4 the following

**Theorem 6.1** *Let $(M, G)$ be a (closed) cohomogeneity one manifold, and let $L \subset G$ be a closed subgroup which acts freely on $M$. Then $N := L \backslash M$ admits metrics of nonnegative Ricci and almost nonnegative sectional curvature which are invariant under the action of $Norm_G(L)/L$, i.e. for every $\varepsilon > 0$, there is a $Norm_G(L)/L$-invariant metric $g_\varepsilon$ on $N$ such that $Ric(N, g_\varepsilon) \geq 0$ and $Sec(N, g_\varepsilon) \cdot diam(N, g_\varepsilon)^2 > -\varepsilon$.*

**Proof.** By the preceding remarks, we may assume that $N$ is obtained by glueing together two disc bundles of the form considered in Corollary 5.4. The metrics $g_\delta$ constructed there extend to a metric on $N$, also denoted by $g_\delta$, since along the boundary both metrics are isometric to $dt^2 + g_Q$.

Thus, we have $Sec(N, g_\delta) \geq -\delta$, $Ric(N, g_\delta) \geq 0$ and $diam(N, g_\delta) = O(\delta^{-1/6})$, so that $Sec(N, g_\delta) \cdot diam(N, g_\delta)^2 \geq -\delta\, O(\delta^{-1/6})^2 = -O(\delta^{2/3})$ which converges to zero. ∎

Let us now investigate under which circumstances the metrics $g_\varepsilon$ from Theorem 6.1 have points of positive Ricci curvature. This part is essentially an adaptation of the arguments in [GZ2] where cohomogeneity one manifolds are treated. First, we obtain the following result.

**Corollary 6.2** *Let $(M, G)$ be a cohomogeneity one manifold with two singular orbits, and let $L \subset G$ act freely on $M$. Further, suppose that at least one of the $L$-quotients of the singular orbits has finite fundamental group and codimension $\geq 2$. Then $N := L \backslash M$ admits $Norm_G(L)/L$-invariant metrics of positive Ricci and almost nonnegative sectional curvature.*

**Proof.** By the proof of Theorem 6.1, $(N, g_\delta)$ is obtained by glueing together two homogeneous disc bundles $\underline{D}_R^\pm \to L \backslash G/K_\pm$, and thus by Corollary 5.4 and our hypothesis, $(N, g_\delta)$ contains points of positive Ricci curvature. Now one uses the deformation results from [Au], [Eh], [We] to obtain the metric with the asserted properties. ∎

Evidently, this result implies that under the hypotheses stated there $N$ has finite fundamental group. Note that if $M$ and hence $N$ fiber over $S^1$, then both $M$ and $N$ have infinite fundamental group. Thus, if $N$ is closed and has finite fundamental group, then $M$ must have two singular orbits. We now wish to obtain further criteria for the finiteness of $\pi_1(N)$.

**Lemma 6.3** *If the singular orbits $G/K_\pm \subset M$ are both exceptional, i.e. have codimension one, then $N$ has infinite fundamental group.*

**Proof.** Suppose that $M$ is given by the group diagram $H \subset \{K_+, K_-\} \subset G$. If $G/K_+ \subset M$ is exceptional, then the unit normal bundle of $G/K_+$ is a double cover of $G/K_+$. It follows that there is a double cover $\pi : \tilde{M} \to M$ where $\tilde{M}$ is the cohomogeneity one manifold with group diagram $H \subset \{K_-, K_-\} \subset G$. Since $\pi$ is $L$-invariant, it induces a double cover $\underline{\pi} : \tilde{N} \to N$ where $\tilde{N} := L \backslash \tilde{M}$.

If $G/K_-$ is exceptional as well, then $\tilde{M}$ is again a cohomogeneity one manifold with two exceptional orbits, so that again $\tilde{N}$ admits a double cover. Thus, inductively we obtain a $2^k$-fold cover of $N$ for all $k$ so that $N$ has infinite fundamental group. ∎

Recall the two $Q$-orthogonal decompositions

$$\mathfrak{g} = \mathfrak{h} \oplus \mathfrak{m}_1^\pm \oplus \mathfrak{m}_2^\pm, \qquad \mathfrak{k}_\pm = \mathfrak{h} \oplus \mathfrak{m}_1^\pm.$$

**Lemma 6.4** *If $N$ has finite fundamental group then $\mathfrak{z}(\mathfrak{g}) \cap \mathcal{H}_e \cap \mathfrak{m}_2^+ \cap \mathfrak{m}_2^- = 0$, where $\mathcal{H}_e := (\mathfrak{h} + \mathfrak{l})_Q^\perp$.*

**Proof.** By Lemma 6.3, $N$ has infinite fundamental group if $\mathfrak{m}_1^+ = \mathfrak{m}_1^- = 0$, so that we may assume that, say, $\mathfrak{m}_1^- \neq 0$. After passing to a double cover of $N$ if necessary, we may also assume that $\mathfrak{m}_1^+ \neq 0$.

If $G$ acts by cohomogeneity one on $M$, then so does its identity component $G_0$, since we assume that there are no exceptional orbits. Hence, we may assume that $G$ is connected. Also, after passing to a finite cover of $G$, we may assume that $G = G_s \times Z$ where $G_s$ is semisimple and $Z$ is a torus whose Lie algebra equals $\mathfrak{z}(\mathfrak{g})$.

Now consider the map $pr_Z \circ \mu : L \times K_+ \times K_- \to Z$, where $\mu$ is multiplication and $pr_Z : G \to Z$ is the projection. This is a homomorphism, hence its image is a compact subgroup $Z_1 \subset Z$. Thus, if we let $G' := G_s \times Z_1 \subset G$, then $G'$ is a normal subgroup which contains $H, K_\pm$ and $L$. It follows that there is a fibration

$$N \longrightarrow G/G' = Z/Z_1 = T^k,$$

and since $\pi_1(N)$ is finite, we must have $k = 0$, i.e. $Z_1 = Z$, i.e. $pr_{\mathfrak{z}(\mathfrak{g})}(\mathfrak{k}_+ + \mathfrak{k}_- + \mathfrak{l}) = \mathfrak{z}(\mathfrak{g})$ which is equivalent to $\mathfrak{z}(\mathfrak{g}) \cap \mathcal{H}_e \cap \mathfrak{m}_2^+ \cap \mathfrak{m}_2^- = 0$. ∎

**Proposition 6.5** *Suppose that $N = L \backslash M$ has finite fundamental group and the $L$-quotients of the singular orbits $\Sigma_\pm := L \backslash G / K_\pm$ both have either infinite fundamental group or codimension one.*

*Then $\dim \mathfrak{m}_1^\pm = 1$ and $\mathfrak{m}_1^+ \cap \mathfrak{m}_1^- = 0$. Moreover, $\dim(\mathfrak{z}(\mathfrak{g}) \cap \mathcal{H}_e) = 2$, $\dim(\mathfrak{z}(\mathfrak{g}) \cap \mathcal{H}_e \cap \mathfrak{m}_2^\pm) = 1$. In particular, $\mathfrak{z}(\mathfrak{g}) \cap \mathcal{H}_e = \mathfrak{t}_+ \oplus \mathfrak{t}_-$ where $\mathfrak{t}_\pm = pr_{\mathfrak{z}(\mathfrak{g}) \cap \mathcal{H}_e}(\mathfrak{k}_\pm)$ is one-dimensional and $pr_{\mathfrak{z}(\mathfrak{g}) \cap \mathcal{H}_e} : \mathfrak{g} \to \mathfrak{z}(\mathfrak{g}) \cap \mathcal{H}_e$ denotes the orthogonal projection.*

**Proof.** By Lemma 6.3 and the hypothesis we conclude that, say, $\Sigma_+$ has infinite fundamental group and thus $\mathfrak{z}(\mathfrak{g}) \cap \mathcal{H}_e \cap \mathfrak{m}_2^+ \neq 0$ by Theorem 2.1. Also, $\mathfrak{z}(\mathfrak{g}) \cap \mathcal{H}_e \cap \mathfrak{m}_2^+ \cap \mathfrak{m}_2^- = 0$ by Lemma 6.4 and thus $\mathfrak{m}_2^+ \not\subset \mathfrak{m}_2^-$ so that $\mathfrak{m}_1^- \not\subset \mathfrak{m}_1^+$. In particular, $\mathfrak{m}_1^- \neq 0$ and hence $\mathfrak{z}(\mathfrak{g}) \cap \mathcal{H}_e \cap \mathfrak{m}_2^- \neq 0$ as well.

Note that $\operatorname{codim}(\Sigma_\pm \subset N) = \dim \mathfrak{m}_1^\pm + 1$. If, say, $\dim \mathfrak{m}_1^+ \geq 2$, then $\operatorname{codim}(\Sigma_+ \subset N) \geq 3$, and since $N \backslash \Sigma_+$ is a disc bundle over $\Sigma_-$, it follows that the inclusion $\Sigma_- \hookrightarrow N$ induces an isomorphism $\pi_1(\Sigma_-) \cong \pi_1(N)$ so that $\Sigma_-$ has finite fundamental group, which is a contradiction.

Thus, we conclude that $\dim \mathfrak{m}_1^+ = 1$ and analogously, $\dim \mathfrak{m}_1^- = 1$. Moreover, since $\mathfrak{m}_1^- \not\subset \mathfrak{m}_1^+$, we have $\mathfrak{m}_1^+ \cap \mathfrak{m}_1^- = 0$. Also, $0 = \mathfrak{z}(\mathfrak{g}) \cap \mathcal{H}_e \cap \mathfrak{m}_2^+ \cap \mathfrak{m}_2^- = \mathfrak{z}(\mathfrak{g}) \cap \mathcal{H}_e \cap (\mathfrak{m}_1^+ \oplus \mathfrak{m}_1^-)^\perp$ implies that $\dim(\mathfrak{z}(\mathfrak{g}) \cap \mathcal{H}_e) \leq 2$. As the intersection of the spaces $\mathfrak{z}(\mathfrak{g}) \cap \mathcal{H}_e \cap \mathfrak{m}_2^\pm \neq 0$ vanishes, they are one-dimensional and $\dim(\mathfrak{z}(\mathfrak{g}) \cap \mathcal{H}_e) = 2$. Finally, the orthogonal complement of $\mathfrak{z}(\mathfrak{g}) \cap \mathcal{H}_e \cap \mathfrak{m}_2^\pm$ in $\mathfrak{z}(\mathfrak{g}) \cap \mathcal{H}_e$ equals $\mathfrak{t}_\pm$, and this finishes the proof. ∎

**Theorem 6.6** *Suppose that $N = L \backslash M$ has finite fundamental group. Then $N$ admits a metric with positive Ricci curvature whose isometry group contains $Norm(L)/L$.*

**Proof.** By Corollary 6.2, we may restrict to the case where the hypotheses of Proposition 6.5 apply. Thus, $\mathfrak{m}_1^\pm \neq 0$, hence the orbits of $G$ and $G_0$ on $M$ coincide, so that we may replace $G$ by any subgroup containing $G_0$. By doing this, we may assume that $G/Norm(L)$ is connected. Also, the principal orbit $G/H$ is connected as $M$ is connected.

The proof follows closely the proofs of Theorems 4.4 and 4.5 in [GZ2]. The first step is to determine a normal subgroup $G' \subset G$ which contains the identity components $H_0$ and $L_0$, and such that $G_0/G' = T^2$. For this, we pass to a finite cover $\pi : G_s \times Z \to G_0$, where $G_s$ is semisimple and $Z$ is a torus with Lie algebra $\mathfrak{z}(\mathfrak{g})$, and we let $\tilde{G}' := G_s \times (pr_Z \circ \mu)(\tilde{H}_0 \times \tilde{L}_0)$ where $\tilde{H}_0, \tilde{L}_0$ are the identity components of the preimages of $H$ and $L$ under $\pi$, where $pr_Z : G_s \times Z \to Z$ is the projection and $\mu$ is the multiplication map. Since $pr_Z \circ \mu$ is a homomorphism, $\tilde{G}'$ is a closed connected subgroup, hence so is its image $G' := \pi(\tilde{G}') \subset G$.

To show that $G' \subset G$ is normal, it suffices to show that its Lie algebra $\mathfrak{g}' \cong \mathfrak{g}_s \oplus pr_{\mathfrak{z}(\mathfrak{g})}(\mathfrak{h} + \mathfrak{l})$ is $Ad_G$-invariant since $G'$ is connected. As $Ad_{G_0}|_{\mathfrak{z}(\mathfrak{g})}$ is trivial, it follows that $Ad_g|_{\mathfrak{z}(\mathfrak{g})}$ depends only on the component of $g \in G$. But each component of $G$ contains elements of $H$ and $Norm(L)$ as $G/H$ and $G/Norm(L)$ are connected. Therefore, $Ad_g|_{\mathfrak{z}(\mathfrak{g})}$ preserves $pr_{\mathfrak{z}(\mathfrak{g})}(\mathfrak{h} + \mathfrak{l})$ and hence $Ad_g(\mathfrak{g}') = \mathfrak{g}'$ for all $g \in G$.

Since $G/K_\pm$ is also connected, the spaces $\mathfrak{t}_\pm = pr_{(\mathfrak{g}')^\perp}(\mathfrak{k}_\pm)$ are also $Ad_G$-invariant, and since $(\mathfrak{g}')^\perp = \mathfrak{z}(\mathfrak{g}) \cap \mathcal{H}_e = \mathfrak{t}_+ \oplus \mathfrak{t}_-$ by Proposition 6.5, it follows that $Ad_G|_{(\mathfrak{g}')^\perp}$ is a subgroup of $\mathbb{Z}_2 \oplus \mathbb{Z}_2$, acting by $\pm Id$ on $\mathfrak{t}_\pm$.

Now one concludes as in the proof of Theorem 4.4 in [GZ2] that $B^3 := L\backslash M/G'$ is an orbifold homeomorphic to $S^3$ with possible orbifold singularities along the singular orbits of the standard $T^2$-action on $S^3 \subset \mathbb{C}^2$. The principal fiber is $(L \cap G')\backslash G'/(H \cap G')$, whereas the fibers over the orbifold singularities are $(L \cap G')\backslash G'/(K_\pm \cap G')$ which are finitely covered by the principal fiber with cyclic groups of deck transformations. Moreover, the action of $G/G'$ on $B^3$ is contained in the extension of the standard $T^2$-action by the $\mathbb{Z}_2 \oplus \mathbb{Z}_2$-action induced by complex conjugation of each complex coordinate of $S^3 \subset \mathbb{C}^2$. Thus, the orbifold metric $g_0$ on $B^3$ induced by the standard metric on $S^3$ is $G/G'$-invariant.

Since $G' \subset G$ is normal, it follows that the tangent distribution $\mathcal{V}_M$ to the $G'$-orbits on $M$ has constant rank and is $G$-invariant, and we can choose a $G$-invariant complementary distribution $\mathcal{H}_M$ (e.g. the orthogonal complement w.r.t. any $G$-invariant metric on $M$). As $L_0 \subset G'$, the distributions $\mathcal{V} := d\pi(\mathcal{V}_M)$ and $\mathcal{H} := d\pi(\mathcal{H}_M)$ on $N$ are well-defined, where $\pi : M \to N$ is the canonical projection. Now we define a $G$-invariant metric $g_M$ on $M$ by fixing the normal homogeneous metric induced by $Q$ on $\mathcal{H}_M \cong T(G'/(H \cap G'))$, declaring $\mathcal{H}_M$ and $\mathcal{V}_M$ to be orthogonal, and such that the canonical orbifold submersion $p_M : (M, g_M) \to (B^3, g_0)$ is Riemannian.

The $G$-invariance of $g_M$ implies that there is a unique metric $g$ on $N$ such that $\pi : (M, g_M) \to (N, g)$ becomes a Riemannian submersion. Note that the $p_M$-fibers of $g_M$ are totally geodesic, hence so are the fibers of the canonical projection $p : N \to B^3$.

We regard $\mathcal{H}$ as an Ehresmann connection on the fiber bundle $p : N \to B^3$ and recall that parallel transport of the fibers along any path in $B^3$ is an isometry iff all fibers are totally geodesic. Thus, if we change the metric on some fiber $F$ to a metric $g'_F$ with the same isometry group as $g|_F$, then parallel translation induces a metric $g'$ on $N$ with totally geodesic fibers such that $g'(\mathcal{H}, \_) = g(\mathcal{H}, \_)$ and $g'|_F = g'_F$. Moreover, the isometries of $(N, g)$ and $(N, g')$ which commute with $p$ coincide.

Note that the $p$-fibers are the biquotients $(L \cap G')\backslash G'/(H \cap G')$ which have finite fundamental group and hence admit $Norm(L)/L$-invariant metrics of positive Ricci curvature by Theorem 2.1, since $\mathfrak{z}(\mathfrak{g}') \cap \mathcal{H}_e = \mathfrak{g}' \cap \mathfrak{z}(\mathfrak{g}) \cap \mathcal{H}_e = 0$. Thus, by the preceding remark, we may deform $g$ to a $Norm(L)/L$-invariant metric $g'$ on $N$ whose $p$-fibers are all totally geodesic and have positive Ricci curvature and such that $(N, g') \to (B^3, g_0)$ is again an orbifold submersion.

This is precisely the setting needed to employ the fiber shrinking from [Na] in order to obtain metrics of positive Ricci curvature on $N$ with the asserted isometries; indeed, for this step the proof of Theorem 4.5 in [GZ2] can be employed verbatim. ∎

**Remark 6.7** In the proof of Theorem 6.6, the deformation of the metric $g$ on $N$ to $g'$ is necessary since the $Q$-induced metric on the biquotient $\Sigma := (L \cap G')\backslash G'/(H \cap G')$ may not have positive Ricci curvature everywhere. If $(\Sigma, g_Q)$ *does* have positive Ricci curvature, then we need not replace $g$ by $g'$, and since $(\Sigma, g_Q)$ has also nonnegative sectional curvature, the shrinking of the fibers yields not only metrics of positive Ricci but also of almost nonnegative sectional curvature by [FY]. In particular, since $(\Sigma, g_Q)$ has positive Ricci curvature if $L = 1$, we obtain the following

**Corollary 6.8** *Let $(M, G)$ be a cohomogeneity one manifold with finite fundamental group. Then $M$ carries $G$-invariant metrics of positive Ricci and almost nonnegative sectional curvature.*

# 7 Brieskorn manifolds and their quotients

Particularly interesting examples of closed cohomogeneity one manifolds are given by the odd-dimensional Brieskorn manifolds (see [Bri], [HM], [Mi]). Given integers $n \geq 2$ and $d \geq 1$, the Brieskorn manifolds $W_d^{2n-1}$ are the $(2n-1)$-dimensional real algebraic submanifolds of $\mathbb{C}^{n+1}$ defined by the equations

$$z_0^d + z_1^2 + \cdots + z_n^2 = 0 \quad \text{and} \quad |z_0|^2 + |z_1|^2 + \cdots + |z_n|^2 = 1.$$

All manifolds $W_d^{2n-1}$ are invariant under the standard linear action of $O(n)$ on the $(z_1, \ldots, z_n)$-coordinates, and the action of $S^1$ via the diagonal matrices of the form $diag(e^{2i\theta}, e^{di\theta}, \ldots, e^{di\theta})$. The resulting action of the product group $G := O(n) \times S^1$ has cohomogeneity one ([HH]) with two nonprincipal orbits of codimensions 2 and $n-1$, respectively.

The topology of the Brieskorn manifolds is fairly well understood. Notice as special cases that $W_1^{2n-1}$ is equivariantly diffeomorphic to $S^{2n-1}$ with the linear action of $G = O(n) \times S^1 \subset U(n)$. Moreover, $W_2^{2n-1}$ is equivariantly diffeomorphic to the Stiefel manifold $V_{n+1,2} = O(n+1)/O(n-1)$ of orthonormal two-frames in $\mathbb{R}^{n+1}$, with the action of $G$ given by the standard inclusion $G \subset O(n+1) \times O(2) = O(n+1) \times Norm_{O(n+1)}O(n-1)/O(n-1)$. Also, for $n = 2$ the manifold $W_d^3$ is diffeomorphic to the lens space $S^3/\mathbb{Z}_d$, whereas for $n \geq 3$ all manifolds $W_d^{2n-1}$ are simply connected.

If now $n \geq 3$ and $d \geq 1$ are odd, then $W_d^{2n-1}$ is a homotopy sphere. Indeed, if $d \equiv \pm 1 \mod 8$ then $W_d^{2n-1}$ is diffeomorphic to the standard $(2n-1)$-sphere, while for $d \equiv \pm 3 \mod 8$, $W_d^{2n-1}$ is diffeomorphic to the *Kervaire sphere* $K^{2n-1}$. The Kervaire spheres are generators of the groups of homotopy spheres which bound parallelizable manifolds, and they all can be realized as Brieskorn manifolds ([Bre], [HM]). The sphere $K^{2n-1}$ is exotic, i.e. homeomorphic but not diffeomorphic to the standard sphere, if $n+1$ is not a power of 2 ([Bro1]). Whether or not $K^{2n-1}$ is diffeomorphic to the standard sphere is unknown if $n+1 = 2^k$ and $k \geq 6$.

On the other hand, if $n = 2m$ is even, then $W_d^{4m-1}$ is a rational homology sphere whose only nontrivial cohomology groups are given by $H^0(W_d^{4m-1}) \cong H^{4m-1}(W_d^{4m-1}) \cong \mathbb{Z}$ and $H^{2m}(W_d^{4m-1}) \cong \mathbb{Z}_d$ ([Bre], p.275).

Let us also recall that the orbit space of a free action of a nontrivial finite cyclic group on a homotopy sphere is called a *homotopy real projective space* if this group has order two, and a *homotopy lens space* otherwise. Notice that homotopy real projective spaces are always homotopy equivalent to standard real projective spaces ([Wa]), whereas the corresponding statement for homotopy lens spaces does in general not hold. For $m \geq 2$, consider now the cyclic subgroup $\mathbb{Z}_m \subset S^1 \subset G$ which hence acts on $W_d^{2n-1}$ by $\alpha(z_0, z_1, \ldots, z_n) := (\alpha^2 z_0, \alpha^d z_1, \ldots \alpha^d z_n)$, where $\alpha$ is an $m$-th root of unity. One verifies that this action is free iff $m$ and $d$ are relatively prime. Moreover, $Norm(\mathbb{Z}_m) = G$ as $\mathbb{Z}_m$ is contained in the center, so that the resulting quotients $W_d^{2n-1}/\mathbb{Z}_m$ are again cohomogeneity one manifolds. The orbit spaces of these free cyclic group actions on the Brieskorn manifolds have been studied in [AB], [Bro2], [Gi1], [Gi2], [Gi3], [Or], [HM], [Lo]. Combining all of this with Corollary 6.8 (which is slightly stronger than the results in [GZ2] and [ST]), we obtain the following

**Corollary 7.1** *The following closed manifolds are of cohomogeneity one and admit invariant metrics of positive Ricci and almost nonnegative sectional curvature.*

1. *All Brieskorn manifolds, hence the Kervaire spheres and infinite families of rational homology spheres in each dimension $4m - 1$ for $m \geq 2$.*

2. *Quotients of the Kervaire spheres by a free $\mathbb{Z}_2$-action; indeed, for any integer $k \geq 1$ this results in at least $4^k$ oriented diffeomorphism types of homotopy $\mathbb{RP}^{4k+1}$.*

3. *Quotients of the Kervaire spheres by free actions of $\mathbb{Z}_m$ for any integer $m \geq 3$; these quotients are homotopy lens spaces which are differentiably distinct from the standard ones in those dimensions in which the Kervaire sphere is exotic.*

The existence of metrics of positive Ricci curvature on the Kervaire spheres was already treated in [Ch] and [BH]. Recently, Sasakian metrics of positive Ricci curvature on the homotopy $\mathbb{RP}^{4k+1}$'s from Corollary 7.1 have been constructed in [BGN]. Also, it was shown in [GZ1] that the four oriented diffeomorphism types of homotopy $\mathbb{RP}^5$'s actually admit invariant metrics of nonnegative sectional curvature.

Let us now leave the cohomogeneity one case and study free actions of positive-dimensional Lie groups on the Brieskorn manifolds.

Suppose that $n = 2m$ is even. In this case, there is a free circle action on $W_d^{4m-1}$ which is given by the action of the circle subgroup $S^1 = Z(U(m)) \subset O(2m)$ where $Z$ denotes the center. Likewise, if we assume that $n = 4m$, then the subgroup $Sp(1) \subset O(4m)$, viewed as scalar multiplication of the unit quaternions on $\mathbb{R}^{4m} \cong \mathbb{H}^m$, acts also freely on $W_d^{8m-1}$.[1] For the quotient manifolds, we obtain the following

**Proposition 7.2** *If $m \geq 2$, then $N := N_d^{4m-2} := W_d^{4m-1}/S^1$ and $\tilde{N} := \tilde{N}_d^{8m-4} := W_d^{8m-1}/Sp(1)$ are simply connected, and their integral cohomology groups and cohomology rings are given by*

$$H^k(N) \cong H^k(\mathbb{CP}^{2m-1}), \quad H^k(\tilde{N}) \cong H^k(\mathbb{HP}^{2m-1}),$$

$$H^*(N) \cong H^*(\tilde{N}) \cong \mathbb{Z}[x,y]/\{x^m - d\,y,\ y^2\},$$

*where $x \in H^2(N)$ ($x \in H^4(\tilde{N})$, respectively) is the Euler class of the principal bundle $S^1 \hookrightarrow W_d^{4m-1} \to N$ ($Sp(1) \hookrightarrow W_d^{8m-1} \to \tilde{N}$, respectively) and $y$ is the generator of $H^{2m}(N)$ ($H^{4m}(\tilde{N})$, respectively). In particular, $N$ and $\tilde{N}$ have the same rational cohomology ring as $\mathbb{CP}^{2m-1}$ and $\mathbb{HP}^{2m-1}$, respectively.*

*Moreover, $N$ and $\tilde{N}$ admit metrics of positive Ricci and almost nonnegative sectional curvature such that the canonical cohomogeneity two actions of $S^1 \times SU(m)$ on $N$ and $S^1 \times Sp(m)$ on $\tilde{N}$ are isometric. Also, w.r.t. these metrics, the fiber bundles $Sp(1)/S^1 = S^2 \hookrightarrow N_d^{8m-2} \to \tilde{N}_d^{8m-4}$ are Riemannian.*

*Furthermore, for $d \geq 3$ the manifolds $N$, $\tilde{N}$ neither admit a Lie group action of cohomogeneity less than two, nor are they diffeomorphic to a biquotient.*

**Proof.** Consider first the quotients $N$. From the homotopy exact sequence of the fibration $S^1 \hookrightarrow W_d^{4m-1} \to N$ it follows that $N$ is simply connected, whence $H^1(N) = H^{4m-3}(N) = 0$. The asserted cohomology ring structure of $H^*(N)$ follows then from the Gysin sequence.

By Corollary 6.2, $N$ admits metrics of positive Ricci and almost nonnegative sectional curvature whose isometry group contains $Norm(S^1)/S^1 = S^1 \times SU(m)$, and it is straightforward to verify that this action has cohomogeneity two.

The corresponding assertions for $\tilde{N}$ follow analogously.

In [Uch] and [Iw], a classification of closed cohomogeneity one manifolds with the rational cohomology ring of $\mathbb{CP}^n$ and $\mathbb{HP}^n$, respectively, was established. Moreover, in [KZ] all biquotients (and homogeneous spaces) whose rational cohomology ring is generated by one element were classified. From these results, the final assertion follows. ∎

From the above description of $W_d^{2n-1}$ for $d = 1, 2$ it follows that $N_1^{4m-2} = \mathbb{CP}^{2m-1}$ and $\tilde{N}_1^{8m-4} = \mathbb{HP}^{2m-1}$ are homogeneous, and the bundle $S^2 \hookrightarrow \mathbb{CP}^{4m-1} \to \mathbb{HP}^{2m-1}$ is the quaternionic twistor fibration. Also, $N_2^{4m-2} = S^1 \backslash O(2m+1)/O(2m-1)$ and $\tilde{N}_2^{8m-4} = Sp(1)\backslash O(4m+1)/O(4m-1)$ are non-homogeneous biquotients ([KZ]). Here, equality signs indicate a diffeomorphism.

Note that for $m = 1$ the action of $Sp(1)$ on $\tilde{N}_d^4$ has cohomogeneity one. In fact, from the group diagram of this action it follows that $\tilde{N}_d^4$ is diffeomorphic to $S^4$, and by the Gysin sequence, the principal $Sp(1)$-bundle $W_d^7 \to \tilde{N}_d^4$ has Euler class $d \in H^4(\tilde{N}_d^4, \mathbb{Z}) \cong \mathbb{Z}$. Thus, up to orientation, these bundles exhaust all non-trivial $S^3$-bundles over $S^4$ as these are determined by their Euler class ([GZ1]). Moreover, again by [GZ1], $W_d^7$ and hence $N_d^6$ even admit metrics of nonnegative sectional curvature, and there is a fibration $S^2 \hookrightarrow N_d^6 \to \tilde{N}_d^4 = S^4$ for all $d$.

---

[1] This action has been pointed out to us by B. Wilking

Université Libre de Bruxelles, CP218, Campus Plaine, Boulevard de Triomphe, B-1050 Bruxelles, Belgium

Mathematisches Institut, Westfälische Wilhelms-Universität Münster, Einsteinstrasse 62, D-48149 Münster, Germany

Email: lschwach@ulb.ac.be, wtusch@math.uni-muenster.de